\documentclass[smallextended,referee,envcountsect,]{svjour3}

\smartqed
\usepackage{comment}
\usepackage{booktabs,framed}
\usepackage{placeins}
\usepackage{algorithm}
\usepackage[misc]{ifsym}
\usepackage{algpseudocode}
\usepackage{subcaption}
\usepackage{array}    
\usepackage{booktabs} 
\usepackage{amsfonts}
\usepackage{graphicx}
\usepackage{tcolorbox,enumerate}
\usepackage{amssymb,latexsym,amsmath,enumitem,amscd,amsthm,graphicx,color}
\usepackage{hyperref}
\allowdisplaybreaks

\newtheorem{assumption}{Assumption}
\begin{document}

\title{Variance-Reduced Stochastic Quasi-Newton SQCQP Method for Nonlinear Programming Problems}

\author{Gopal Tiwari, Md Abu Talhamainuddin Ansary}

\institute{Gopal Tiwari and  Md Abu Talhamainuddin Ansary \at
Department of Mathematics\\
IIT Jodhpur, Jodhpur, India 342030\\
Email:\\
G. Tiwari: gopalpersonalid2003@gmail.com\\
M. A. T, Ansary: md.abutalha2009@gmail.com (\Letter)
 }


\maketitle


\begin{abstract}
In this paper, we propose a variance-reduced stochastic quasi-Newton sequential quadratic constrained quadratic programming method (VR-QN-SQCQP) for finite-sum nonlinear optimization with equality and inequality constraints. Exact nonlinear-constraint Hessians are replaced entirely by safeguarded symmetric rank-one (SR1) curvature models built from first-order constraint gradients, while PAGE estimates the objective gradient and an exact \(L_1\) merit function provides globalization. Under explicit regularity, boundedness, relative QCQP-accuracy, and step-length assumptions, the analysis establishes almost-sure subsequential convergence to Karush--Kuhn--Tucker (KKT) points and, under stronger uniform step-length and error-control conditions, an expected \(O(\varepsilon^{-2})\) outer-iteration bound for the squared KKT residual. A conditional local stability result also compares learned-curvature directions with their exact-Hessian counterparts. Numerical experiments compare VR-QN-SQCQP against Exact-SQCQP-PAGE on Sparse PCA using 253,680 CDC BRFSS 2015 observations and on a dense nonlinear benchmark. Over ten seeds, both methods converge in all Sparse PCA runs, while on the dense benchmark the proposed method converges in 10/10 runs versus 9/10 for the exact-Hessian variant. The proposed method performs zero exact constraint-Hessian evaluations and reduces median curvature-construction time on the dense benchmark from 0.1230 s to 0.0111 s. Overall runtime benefits remain problem-dependent.
\end{abstract}

\section{Introduction}
\label{sec:introduction}

Stochastic nonlinear optimization with deterministic equality and inequality constraints arises in statistical learning, scientific computing, and engineering design. In large-scale applications, the objective is often expressed as an expectation or finite sum, making repeated full-gradient evaluations computationally expensive. Stochastic approximation and variance-reduction methods provide effective alternatives for reducing this cost \cite{robbinsmonro1951,nemirovski2009robust,ghadimiLan2013stochastic,bottou2018optimization,duchi2011adaptive}. However, nonlinear constraints introduce additional difficulties involving feasibility, multiplier estimation, and constrained stationarity, motivating methods that combine efficient gradient estimation with reliable constraint handling \cite{bertsekas1982constrained,bonnans2006numerical}.

Sequential quadratic programming (SQP) is a well-established framework for nonlinear constrained optimization, using local quadratic objective models and linearized constraints to generate search directions \cite{powell1978sqp,boggs1995sqp}. Its practical effectiveness depends on suitable globalization strategies, including trust-region methods, filter techniques, and exact penalty functions \cite{conn2000trust,fletcherleyffer2002filter,dipillogrippo1985exact}. Sequential quadratically constrained quadratic programming (SQCQP) extends this approach by retaining quadratic models of nonlinear constraints within each subproblem. Global and local convergence properties of SQCQP methods have been investigated by Fukushima, Luo, and Tseng, and by Solodov \cite{fukushima2003sqcqp,solodov2004sqcqp}. This structure provides a means of incorporating nonlinear constraint curvature directly into the computation of search directions.

An issue with solving SQCQP problems numerically is that calculating the hessian of the nonlinear constraints can be quite costly. However, to overcome this difficulty we can use quasi-Newton methods that approximate the curvature information using the difference between gradients \cite{dennismore1977qn,nocedal2006numerical}. This method has been shown to be useful for large scale SQP implementations \cite{gill2005snopt}. However, approximating the objective or Lagrangian hessian is not the same as approximating each nonlinear constraint's curvature individually. It is then interesting to see if we can replace all the exact constraint Hessians with learned approximations, and still gain the computational benefits whilst keeping the useful local properties of SQCQP.In addition, estimating gradients of objectives is a complementing problem. Stochastic and adaptive sampling schemes have been proposed to optimize finite sums \cite{byrd2012sample}. Variance reduction techniques such as SVRG \cite{johnsonzhang2013svrg}, SAGA \cite{defazio2014saga}, SARAH \cite{nguyen2017sarah} and SPIDER \cite{fang2018spider} have also been studied. The Probabilistic Gradient Estimator (PAGE) algorithm includes random full gradient refreshes and recursive minibatch corrections with complexity guarantees for non-convex optimization \cite{li2021page}.In this work we make use of PAGE as a variance reduction technique. The combination with constrained globalization is non-trivial because the search direction is computed using the gradient approximation; however, the merit function is evaluated using the true objective and constraint function values. This means descent properties of the model gradient do not necessarily translate to the true merit function decreasing sufficiently. Recent work within the probabilistic trust region and stochastic line search frameworks highlights the importance of having explicit control of uncertainty in the model \cite{bandeira2014probabilistic,paquette2020line}.

Recent work has established important foundations for stochastic sequential quadratic optimization under nonlinear constraints. This includes stochastic SQO methods for equality constraints, adaptive methods employing exact augmented Lagrangian merit functions, active-set methods for inequality constraints, approaches for handling rank-deficient constraints, and convergence guarantees with high probability \cite{berahas2021sqo,na2023adaptive,na2023active,berahas2024rankdeficient,curtis2024sqo,berahas2025highprob}. In related work, the authors have developed an adaptive variance-reduced stochastic sequential (quadratic constrained) quadratic programming framework, combining SPIDER-based gradient estimation with a relative constraint-model-curvature indicator to adaptively select between SQP and SQCQP subproblems \cite{TiwariAnsary2026AVRSSQCQP}. While this framework adapts the subproblem formulation according to the estimated importance of constraint curvature, the present work investigates a complementary question: whether all evaluations of nonlinear-constraint Hessians can be eliminated while retaining quadratic constraint models within the SQCQP framework. Specifically, separate quasi-Newton models are maintained for the curvature of individual constraints, and the resulting method is analyzed to establish convergence guarantees and assess its local behaviour relative to an SQCQP method employing exact constraint Hessians.

A variance reduced quasi-Newton SQCQP method, referred to as VR-QN-SQCQP, is developed. Separate symmetric rank one (SR1) matrices are used to approximate the curvature of the non-linear constraints using deterministic gradient difference. PAGE is used to approximate the gradient of the objective function and a safeguarded SR1 model is used to approximate curvature. Trial directions are generated by solving regularized SQCQP subproblems. Globalization is achieved using an exact \(L_1\) merit function with an Armijo backtracking strategy. Three main contributions are made. Firstly, we develop a general framework to replace all the exact non-linear constraint Hessians with separate SR1 models. Secondly, a conditional convergence analysis is presented which takes into account the error incurred by using gradient approximations, and the approximation involved in solving subproblems. We provide convergence guarantees for subsequences of iterates to satisfy KKT conditions under the given assumptions. Additionally, stronger results are derived for residual convergence and complexity given some extra conditions. Finally, a local stability result is provided under suitable curvature consistency and regularity assumptions. Numerical experiments are performed comparing with Exact-SQCQP-PAGE which maintains the gradient estimation and globalization strategy but keeps exact Hessians. The methods converge for all runs of Sparse PCA for all 10 seeds. For the dense non-linear test problem, VR-QN-SQCQP converges for all 10 runs whereas Exact-SQCQP-PAGE only converges for 9 runs. Additionally, the median time taken to build the constraint curvature is reduced significantly from $0.1230$ s to $0.0111$ s on the dense test problem. However, the benefit to runtimes is dependent on the problem.

The remainder of the paper is organized as follows. Section~\ref{sec:preliminaries} introduces the problem formulation and notation. Section~\ref{sec:algorithm} presents the proposed VR-QN-SQCQP method and its globalization strategy. Section~\ref{sec:convergence} develops the convergence analysis and local curvature-consistency results. Section~\ref{sec:experiments} reports the numerical experiments, and Section~\ref{sec:conclusion} concludes the paper.

\section{Preliminaries and Problem Formulation}
\label{sec:preliminaries}

Consider the stochastic nonlinear constrained optimization problem
\begin{equation}\label{eq:problem}
\begin{aligned}
\min_{x\in\mathbb{R}^{n}}\quad
& f(x)=\mathbb{E}_{\xi}[F(x;\xi)]\\
\mathrm{s.\ t.}\quad
& h_i(x)=0,\qquad i=1,\ldots,m,\\
& c_j(x)\leq0,\qquad j=1,\ldots,p,
\end{aligned}\end{equation}
where $F(x;\xi)$ denotes a stochastic objective function associated with a random variable $\xi$, $h_i:\mathbb{R}^{n}\rightarrow\mathbb{R}$, $i=1,\ldots,m$, represent nonlinear equality constraints, and $c_j:\mathbb{R}^{n}\rightarrow\mathbb{R}$, $j=1,\ldots,p$, represent nonlinear inequality constraints. The computational framework developed in this paper focuses on the finite-sum specialization $f(x)=\frac{1}{N}\sum_{\ell=1}^{N}F_{\ell}(x),$
where $N$ denotes the number of samples and each $F_{\ell}:\mathbb{R}^{n}\rightarrow\mathbb{R}$ is differentiable. The equality and inequality constraint vectors are defined by
$h(x)=
\begin{bmatrix}
h_1(x) & \cdots & h_m(x)\end{bmatrix}^{\top},
c(x)=
\begin{bmatrix}
c_1(x) & \cdots & c_p(x)
\end{bmatrix}^{\top},
$
with corresponding Jacobian matrices
$
J_h(x)=
\begin{bmatrix}
\nabla h_1(x) & \cdots & \nabla h_m(x)
\end{bmatrix}^{\top},
J_c(x)=
\begin{bmatrix}
\nabla c_1(x) & \cdots & \nabla c_p(x)
\end{bmatrix}^{\top}
$ respectively.

Throughout the paper, $\|\cdot\|$ denotes the Euclidean norm for vectors and the spectral norm for matrices, unless otherwise specified. The Frobenius norm is denoted by
$\|\cdot\|_F$, while $[v]^+$ represents the componentwise positive part of a vector $v$. The space of real symmetric $n\times n$ matrices is denoted by $\mathbb{S}^{n}$.

\begin{definition}[Karush--Kuhn--Tucker Conditions]
\label{def:kkt_conditions}
A feasible point $x^*$ of~\eqref{eq:problem} is said to satisfy the Karush--Kuhn--Tucker (KKT) conditions if there exist multipliers
$\lambda^*\in\mathbb{R}^{m}$ and
$\mu^*\in\mathbb{R}^{p}_{+}$ such that
\begin{equation}\label{eq:kkt}
\begin{aligned}
\nabla f(x^*)
+J_h(x^*)^{\top}\lambda^*
+J_c(x^*)^{\top}\mu^*&=0,\\
h(x^*)&=0,\\
c(x^*)&\leq0,\\
\mu^*&\geq0,\\
\mu_j^*c_j(x^*)&=0,
\qquad j=1,\ldots,p.
\end{aligned}
\end{equation}
Under a suitable constraint qualification, these conditions are necessary for local optimality.
\end{definition}

\begin{definition}[Linear Independence Constraint Qualification]\label{def:licq} The active inequality index set at a feasible point $x$ is defined by $\mathcal{A}(x)=\left\{j\in\{1,\ldots,p\}:c_j(x)=0\right\}.$ The linear independence constraint qualification (LICQ) is said to hold at $x$ if the vectors $\left\{\nabla h_i(x)\right\}_{i=1}^{m}\cup\left\{\nabla c_j(x)\right\}_{j\in\mathcal{A}(x)}$ are linearly independent.
\end{definition}

To quantify approximate stationarity
and feasibility, the KKT residual mapping
is introduced as
\begin{equation}\label{eq:kkt_residual}
\mathcal{R}(x,\lambda,\mu)=
\begin{bmatrix}
\nabla f(x)
+J_h(x)^{\top}\lambda
+J_c(x)^{\top}\mu\\
h(x)\\
[c(x)]^+\\
\mu\odot c(x)
\end{bmatrix},
\end{equation}
where $\odot$ denotes component-wise multiplication. For inequality multipliers satisfying $\mu\geq0$, the corresponding KKT residual norm is defined by
$\mathcal{E}(x,\lambda,\mu)=\|\mathcal{R}(x,\lambda,\mu)\|_2.$

For a twice continuously differentiable constraint function $r$, the second-order Taylor expansion at $x^k$ is 
$r(x^k+d) =r(x^k)+\nabla r(x^k)^{\top}d +\frac12d^{\top}\nabla^2r(x^k)d +o(\|d\|^2)$. To avoid exact constraint-Hessian evaluations, symmetric quasi-Newton matrices
$C_{h_i}^k,C_{c_j}^k\in\mathbb{S}^{n}$ are introduced. The resulting quadratic constraint models are
\begin{equation}\label{eq:constraint_models}
\begin{aligned}
\widehat h_i^k(d)
={}&h_i(x^k)+\nabla h_i(x^k)^{\top}d
+\frac12d^{\top}C_{h_i}^kd,\\
\widehat c_j^k(d)
={}&c_j(x^k)+\nabla c_j(x^k)^{\top}d
+\frac12d^{\top}C_{c_j}^kd.
\end{aligned}
\end{equation}

For globalization, define the total constraint violation  $V(x)=\|h(x)\|_1+\|[c(x)]^+\|_1,$ and for any penalty parameter $\rho>0$, the exact $L_1$ merit function is
\begin{equation}\label{eq:merit}
\Phi(x;\rho)=f(x)+\rho V(x).\end{equation}

Let $g^k$ denote the PAGE stochastic objective-gradient estimator and define $e^k=g^k-\nabla f(x^k)$. For any direction $d$ and penalty $\rho>0$, the exact and stochastic merit directional quantities are, respectively,
$\Phi'(x^k;d,\rho)=\nabla f(x^k)^\top d+\rho V'(x^k;d)$ and
$D_{\Phi,k}(d;\rho)=(g^k)^\top d+\rho V'(x^k;d).$

Consequently,
\begin{equation}\label{eq:merit_direction_error}
D_{\Phi,k}(d;\rho)-\Phi'(x^k;d,\rho)=(e^k)^\top d. \end{equation}

In the proposed method, $\rho_k$ denotes the stored penalty, whereas $\varrho_k$ denotes the trial penalty used in a successful Armijo test. The construction and commitment of
these penalties are described in Section~\ref{sec:algorithm}.

\section{Proposed Variance-Reduced Stochastic
Quasi-Newton SQCQP Method}
\label{sec:algorithm}

A variance-reduced stochastic quasi-Newton sequential quadratic constrained quadratic programming (VR-QN-SQCQP) method is developed for problem~\eqref{eq:problem}. 

The objective gradient is estimated using PAGE\cite{li2021page}. Let $g^k\in\mathbb R^n$ denote the stochastic gradient estimator used to construct the accepted QCQP model at iteration $k$, with initialization
$g^0=\nabla f(x^0)$. Let $\mathcal B_{\mathrm s}^k$ be a uniformly sampled mini-batch of size $b_{\mathrm s}$.
Following an accepted nonzero step, the PAGE estimator is updated according to
\begin{equation}\label{eq:page_update}
g^{k+1}=
\begin{cases}
\nabla f(x^{k+1}),
&\text{with probability }p_{\mathrm{ref}},\\[2mm]
g^k+\dfrac{1}{b_{\mathrm s}}
\displaystyle\sum_{\ell\in\mathcal B_{\mathrm s}^k}
\left[
\nabla F_\ell(x^{k+1})
-\nabla F_\ell(x^k)
\right],
&\text{with probability }1-p_{\mathrm{ref}}.
\end{cases}
\end{equation}

To approximate objective and constraint curvature, safeguarded SR1 updates are employed \cite{nocedal2006numerical}. Let $G^k\in\mathbb S^n$ denote the stored objective curvature approximation, initialized by $G^0=I$. After an accepted nonzero step, define $s^k=x^{k+1}-x^k,$ and when a PAGE mini-batch correction is performed, the objective secant quantities are $ y_f^k=g^{k+1}-g^k,  v_f^k=y_f^k-G^ks^k.$ The objective curvature is updated by
\begin{equation}
\label{eq:objective_sr1}
G^{k+1}=
\begin{cases}
G^k+
\dfrac{v_f^k(v_f^k)^\top}
{(v_f^k)^\top s^k},
&\text{if safeguards hold},\\[3mm]
G^k,
&\text{otherwise}.
\end{cases}
\end{equation}
An update is accepted only if $|(v_f^k)^\top s^k|\geq\kappa_G\|v_f^k\|\|s^k\|$ and the candidate matrix satisfies the prescribed boundedness and secant safeguards.
Objective SR1 updates are skipped during PAGE full-gradient refreshes to avoid forming secants from inconsistent stochastic
gradient information. The stored matrix $G^k$ is permitted to be indefinite, thereby retaining potentially useful nonconvex objective curvature.

For each nonlinear constraint $r\in\{h_1,\ldots,h_m,c_1,\ldots,c_p\}$, an independent curvature matrix $C_r^k\in\mathbb S^n$ is maintained, with initialization $C_r^0=0$. Using deterministic constraint gradients, define
$y_r^k=\nabla r(x^{k+1})-\nabla r(x^k), v_r^k=y_r^k-C_r^ks^k.$
The corresponding SR1 update is
\begin{equation}\label{eq:constraint_sr1}
C_r^{k+1}=
\begin{cases}
C_r^k+
\dfrac{v_r^k(v_r^k)^\top}
{(v_r^k)^\top s^k},
&\text{if safeguards hold},\\[3mm]
C_r^k,
&\text{otherwise}.
\end{cases}
\end{equation}
The update is accepted only when
$|(v_r^k)^\top s^k|\geq\kappa_C\|v_r^k\|\|s^k\|, \
\|C_{r,\mathrm{trial}}^{k+1}\|_F
\leq M_r.
$
The Frobenius-norm safeguard requires $O(n^2)$ arithmetic operations for each dense constraint-SR1 update and ensures
$\|C_r^k\|_2\leq\|C_r^k\|_F\leq M_r,$ provided the initialization satisfies the same bound. At iteration $k$, the stochastic gradient estimator and quasi-Newton curvature models
are incorporated into a regularized SQCQP subproblem. Using the quadratic constraint models in~\eqref{eq:constraint_models}, a trial direction is obtained from
\begin{equation}\label{eq:qn_sqcqp_subproblem}
\begin{aligned}
\min_{d\in\mathbb R^n}\quad&
(g^k)^\top d+
\frac12d^\top(G^k+\zeta_k I)d\\
\mathrm{s.~t.~}\quad&
\widehat h_i^k(d)=0,
\quad i=1,\ldots,m,\\
&
\widehat c_j^k(d)\leq0,
\quad j=1,\ldots,p.
\end{aligned}
\end{equation}
Here $\zeta_k\geq0$ is a solve-time regularization parameter. Initially, $\zeta_k=0$ is attempted and may be increased during numerical recovery. The regularized matrix $G^k+\zeta_k I$ is used only in the current QCQP; the stored approximation remains $G^k$.

Since the curvature matrices may be indefinite, the QCQP is generally nonconvex. Consequently, a numerically acceptable approximate local solution is sought rather than assuming global subproblem optimality. The numerical implementation employs an SLSQP-based solver \cite{kraft1988sqp}.
Let $d^k$ denote the finally selected trial direction and
$(\widehat\lambda^k,\widehat\mu^k)$ its associated recovered QCQP multipliers, where $\widehat\lambda^k\in\mathbb R^m$
and $\widehat\mu^k\in\mathbb R_+^p$.

The solver evaluates local feasibility, stationarity, dual-feasibility, and complementarity residuals. However
, numerical acceptance of an ordinary nonzero direction does not imply the asymptotic relative QCQP accuracy condition used
in Section~\ref{sec:convergence}. That stronger condition is imposed separately in the theoretical analysis.

Globalization is performed using the exact $L_1$ merit function
$\Phi(x;\rho)$ defined in~\eqref{eq:merit}. An important distinction is made between the stored penalty $\rho_k$ entering iteration $k$ and the trial penalty $\varrho_k$ used to evaluate a candidate direction.

For each QCQP candidate successfully computed, let  $M_k=\max\left\{\|\widehat\lambda^k\|_\infty,\|\widehat\mu^k\|_\infty\right\}.$
The candidate penalty is initialized from the stored value $\rho_k$ and increased geometrically, when necessary,
until $\varrho_k\geq(1+\theta_\rho)M_k,\qquad\varrho_k\geq\rho_k,$ where $\theta_\rho>0$ is the penalty margin. The candidate value is used for predicted descent and Armijo evaluation, but is not immediately assigned to the stored penalty.

Consequently, rejected candidate directions do not permanently commit their trial penalty increases. Recovery may produce different trial penalties for different QCQP candidates within
the same outer iteration.

For a candidate direction $d$, the stochastic merit directional quantity is
$D_{\Phi,k}(d;\varrho_k)=(g^k)^\top d+\varrho_k V'(x^k;d).$
Every ordinary accepted nonzero direction is required to satisfy the sufficient predicted-descent condition
\begin{equation}\label{eq:merit_descent}
D_{\Phi,k}(d^k;\varrho_k)\leq-\gamma_\Phi\|d^k\|^2, \end{equation} where $\gamma_\Phi>0$.

If a trial direction fails the required numerical QCQP validation or sufficient merit descent, the subproblem may be resolved using a cold start, increased regularization, or a refreshed PAGE estimator at the current iterate. Regularization and gradient recovery are numerical safeguards rather than unconditional guarantees of successful subproblem solution.

Following successful direction and descent validation, Armijo backtracking \cite{armijo1966minimization} seeks
$\alpha_k\in\{1,\beta,\beta^2,\ldots\}$
such that
\[\begin{aligned}
\Phi(x^k+\alpha_kd^k;\varrho_k)\leq{}&\Phi(x^k;\varrho_k)+\sigma\alpha_k D_{\Phi,k}(d^k;\varrho_k),\end{aligned}\]
where $\sigma,\beta\in(0,1)$. The accepted primal iterate is then $x^{k+1}=x^k+\alpha_kd^k.$

The stochastic gradient error $e^k=g^k-\nabla f(x^k)$ introduces a distinction between predicted and exact merit descent. From~\eqref{eq:merit_direction_error} and~\eqref{eq:merit_descent},
\[\begin{aligned}\Phi'(x^k;d^k,\varrho_k)&=D_{\Phi,k}(d^k;\varrho_k)-(e^k)^\top d^k\\&\leq-\gamma_\Phi\|d^k\|^2
+\|e^k\|\|d^k\|.\end{aligned}\]
Therefore, predicted stochastic descent does not independently guarantee descent of the exact merit function. The convergence analysis explicitly accounts for stochastic-gradient errors, accepted step lengths, and approximate QCQP solution accuracy.

After an ordinary accepted nonzero step, the QCQP multipliers are transferred using the damped updates
\begin{equation}\label{eq:multiplier_update}\begin{aligned}
\lambda^{k+1}&=\lambda^k+\alpha_k(\widehat\lambda^k-\lambda^k),\\\mu^{k+1}&=\left[\mu^k+\alpha_k(\widehat\mu^k-\mu^k)\right]^+.
\end{aligned}\end{equation}
The update becomes a full multiplier handoff when $\alpha_k=1$, while shorter accepted steps produce damped transfers.

After the multiplier handoff, the stored penalty for the next iteration is computed from the retained pair $(\lambda^{k+1},\mu^{k+1})$. Specifically, starting from $\rho_k$, the same geometric penalty-increase
rule is applied until $\rho_{k+1}\geq(1+\theta_\rho)\max\left\{
\|\lambda^{k+1}\|_\infty,\|\mu^{k+1}\|_\infty
\right\}.$ Thus, $\rho_{k+1}\geq\rho_k$, but generally ${\varrho_k\neq\rho_{k+1}.}$ The successful trial penalty and the committed penalty are therefore treated as distinct sequences in the convergence analysis.

A separate safeguard is applied when the QCQP solver returns a negligible direction. Such a direction is accepted as a terminal case only after additional validation of the KKT conditions at the unchanged outer iterate, using the recovered local multipliers. For this special case, the primal iterate remains unchanged, and a full local multiplier handoff may be performed. Following an ordinary accepted step,
the constraint SR1 matrices are updated using deterministic gradient differences. The PAGE estimator is subsequently advanced, after which objective SR1 is applied only for a mini-batch recursive correction.

The complete procedure is summarized
in Algorithm~\ref{alg:vr_qn_sqcqp}.

\begin{algorithm}[H]
\caption{Variance-Reduced Stochastic Quasi-Newton SQCQP
(VR-QN-SQCQP)}
\label{alg:vr_qn_sqcqp}
\begin{algorithmic}[1]

\Require
$x^0\in\mathbb R^n$,
$\lambda^0\in\mathbb R^m$,
$\mu^0\in\mathbb R_+^p$,
$\rho_0>0$;
$G^0=I$, $C_r^0=0$ for $r\in\mathcal R$;
$p_{\mathrm{ref}},b_{\mathrm s}$;
$\kappa_G,\kappa_C,\{M_r\}_{r\in\mathcal R}$;
$\gamma_\Phi>0$, $\sigma,\beta\in(0,1)$,
$\theta_\rho>0$, $\tau_\rho>1$,
$\varepsilon_{\mathrm{KKT}}>0$.

\State $g^0\gets\nabla f(x^0)$

\For{$k=0,1,\ldots$}

    \If{$\mathcal E(x^k,\lambda^k,\mu^k)
    \leq\varepsilon_{\mathrm{KKT}}$}
        \State \textbf{break}
    \EndIf

    \State Obtain $(d^k,\widehat\lambda^k,\widehat\mu^k)$
    from a local solution of~\eqref{eq:qn_sqcqp_subproblem}

    \State $M_k\gets
    \max\{\|\widehat\lambda^k\|_\infty,
    \|\widehat\mu^k\|_\infty\}$

    \State $\ell_k\gets
    \min\{\ell\in\mathbb N_0:
    \rho_k\tau_\rho^\ell\geq
    (1+\theta_\rho)M_k\}$

    \State $\varrho_k\gets\rho_k\tau_\rho^{\ell_k}$

    \State $D_k\gets D_{\Phi,k}(d^k;\varrho_k)$
    \Comment{Require~\eqref{eq:merit_descent}}

    \State $j_k\gets
    \min\left\{j\in\mathbb N_0:
    \begin{aligned}
    &\Phi(x^k+\beta^j d^k;\varrho_k)\\[-1mm]
    &\leq\Phi(x^k;\varrho_k)
    +\sigma\beta^jD_k
    \end{aligned}
    \right\}$

    \State $\alpha_k\gets\beta^{j_k}$

    \State $s^k\gets\alpha_kd^k$

    \State $x^{k+1}\gets x^k+s^k$

    \ForAll{$r\in\mathcal R$}

        \State $y_r^k\gets
        \nabla r(x^{k+1})-\nabla r(x^k)$

        \State $v_r^k\gets y_r^k-C_r^ks^k$

        \State $C_r^{k+1}\gets
        \operatorname{SR1}(C_r^k,s^k,v_r^k)$
        \Comment{Update~\eqref{eq:constraint_sr1}}

    \EndFor

    \State Draw $u_k\sim\operatorname{Unif}[0,1]$

    \If{$u_k\leq p_{\mathrm{ref}}$}

    \State $g^{k+1}\gets\nabla f(x^{k+1}),$ $G^{k+1}\gets G^k$

    \Else

        \State Sample $\mathcal B_{\mathrm s}^k$,
        $|\mathcal B_{\mathrm s}^k|=b_{\mathrm s}$

        \State $g^{k+1}\gets
        g^k+\dfrac{1}{b_{\mathrm s}}
        \displaystyle\sum_{\ell\in\mathcal B_{\mathrm s}^k}
        \left[
        \nabla F_\ell(x^{k+1})-\nabla F_\ell(x^k)
        \right]$
        \Comment{PAGE update~\eqref{eq:page_update}}

        \State $y_f^k\gets g^{k+1}-g^k,$ $v_f^k\gets y_f^k-G^ks^k$

        \State $G^{k+1}\gets
        \operatorname{SR1}(G^k,s^k,v_f^k)$
        \Comment{Update~\eqref{eq:objective_sr1}}

    \EndIf

    \State $\lambda^{k+1}\gets
    \lambda^k+\alpha_k(\widehat\lambda^k-\lambda^k),$ $\mu^{k+1}\gets
    \left[
    \mu^k+\alpha_k(\widehat\mu^k-\mu^k)
    \right]^+$
    \Comment{Update~\eqref{eq:multiplier_update}}

    \State $q_k\gets
    \min\left\{q\in\mathbb N_0:
    \begin{aligned}
    \rho_k\tau_\rho^q\geq{}&(1+\theta_\rho)\\[-1mm]
    &{}\cdot\max\{\|\lambda^{k+1}\|_\infty,
    \|\mu^{k+1}\|_\infty\}
    \end{aligned}
    \right\}$

    \State $\rho_{k+1}\gets\rho_k\tau_\rho^{q_k}$

\EndFor

\end{algorithmic}
\end{algorithm}

\section{Convergence Analysis}\label{sec:convergence}

We first formulate the assumptions needed for the global analysis.

\begin{assumption}[Global regularity and globalization]
\label{ass:global}
The following conditions hold along an infinite sequence of accepted nonterminal iterations.

\begin{enumerate}
\item[(G1)]
The iterates belong to a compact set $\mathcal K\subset\mathbb R^n$. The objective and constraint functions are twice
continuously differentiable on an open neighborhood
containing $\mathcal K$ and the relevant trial segments.
Their gradients are Lipschitz continuous there. Moreover, $f$ is bounded below by $f_{\min}$. The component gradients satisfy
\[\|\nabla F_\ell(x)-\nabla F_\ell(y)\|\leq L_F\|x-y\|.\]

\item[(G2)]
The objective and constraint curvature models are symmetric and uniformly bounded:
\[\|G^k\|\leq M_G,\qquad\|C_r^k\|\leq M_{C,r},\]
for every constraint $r$. The accepted regularization parameters satisfy $0\leq\zeta_k\leq\overline\zeta$.
The local multipliers satisfy
\[\|\widehat\lambda^k\|_\infty\leq\Lambda,\qquad 0\leq\widehat\mu_j^k\leq\Lambda.\]

\item[(G3)]
Each accepted nonterminal iteration satisfies the sufficient-descent condition and the exact-merit Armijo inequality stated above, with a fixed $\gamma_\Phi>0$.

\item[(G4)]
The accepted trial penalties satisfy
\[0<\underline\varrho\leq\varrho_k\leq\overline\varrho<\infty, \]
and their total variation is bounded by a deterministic constant:
\[
\sum_{k=0}^{\infty}|\varrho_{k+1}-\varrho_k| \leq B_{\mathrm{tr}}<\infty.\]

\item[(G5)]
The PAGE estimator uses conditionally independent sampling with small-batch size $b_s$ and full-gradient refresh probability $p_{\mathrm{ref}}\in(0,1]$. The refresh branch evaluates the exact finite-sum gradient, and any additional forced refresh also resets the estimator to that gradient.

The PAGE error sequence, defined using the estimator actually supplied to the accepted QCQP model, satisfies the conditional recursion
\begin{equation}\label{eq:page_error_recursion}\begin{aligned}
\mathbb E[\|e^{k+1}\|^2\mid\mathcal F_k]\leq{}&
(1-p_{\mathrm{ref}})\|e^k\|^2\\&+\frac{(1-p_{\mathrm{ref}})L_F^2}{b_s}\|x^{k+1}-x^k\|^2.
\end{aligned}\end{equation}
Here $\mathcal F_k$ contains the information used to determine the accepted displacement.
\end{enumerate}
\end{assumption}

\subsection{Inexact QCQP Optimality and Merit Descent}

Define the residual of the local QCQP stationarity
equation by
\[
\begin{aligned}
u^k={}&g^k+G_\zeta^kd^k+\sum_{i=1}^{m}\widehat\lambda_i^k
\bigl(\nabla h_i(x^k)+C_{h_i}^kd^k\bigr)+\sum_{j=1}^{p}\widehat\mu_j^k
\bigl(\nabla c_j(x^k)+C_{c_j}^kd^k\bigr).
\end{aligned}
\]

Let $\delta_k$ be a nonnegative number bounding
the following local QCQP KKT residual components:
\begin{equation}
\label{eq:qcqp_accuracy}
\begin{aligned}
\|u^k\|&\leq\delta_k, \ \\|\widehat h^k(d^k)\|_\infty
&\leq\delta_k,\\
\|[\widehat c^k(d^k)]^+\|_\infty
&\leq\delta_k,\\
\|[-\widehat\mu^k]^+\|_\infty
&\leq\delta_k,\\
\|\widehat\mu^k\odot
\widehat c^k(d^k)\|_\infty
&\leq\delta_k.
\end{aligned}
\end{equation}

The following stronger accuracy condition is used only where explicitly invoked.

\begin{assumption}[Relative QCQP accuracy] \label{ass:relative_qcqp}
There exists $\kappa_\delta>0$ such that
\begin{equation}\label{eq:relative_qcqp_accuracy}
\delta_k\leq\kappa_\delta\frac{\|d^k\|^2}{1+\|d^k\|}
\end{equation}
for all sufficiently large accepted nonterminal iterations.
\end{assumption}
For later use, define $ M_C^\Sigma=\sum_{i=1}^{m}M_{C,h_i}+\sum_{j=1}^{p}M_{C,c_j},$
and
$\ell_k=\max\left\{\|\widehat\lambda^k\|_\infty, \|\widehat\mu^k\|_\infty\right\}.$

\begin{proposition}[A sufficient predicted-descent condition]
\label{prop:inexact_qcqp_descent}
Suppose that (G2) holds and the residual bounds \eqref{eq:qcqp_accuracy} are satisfied.

Then the stochastic merit directional quantity obeys
\[ \begin{aligned}
D_{\Phi,k}(d^k;\varrho_k)\leq{}&-\left[\zeta_k-M_G \frac{\varrho_k+\ell_k}{2}M_C^\Sigma\right]\|d^k\|^2-(\varrho_k-\ell_k)V(x^k)+\delta_k\left(\|d^k\|+m\ell_k+p
+\varrho_k(m+p)\right).\end{aligned} \]

In particular, if $\varrho_k\geq\ell_k$, Assumption~\ref{ass:relative_qcqp} holds, and
\[\begin{aligned}\zeta_k\geq{}&M_G+\frac{\varrho_k+\ell_k}{2}
M_C^\Sigma+\gamma_\Phi+\kappa_\delta\max\{1,A_k\},\end{aligned}\]
where $A_k=m\ell_k+p+\varrho_k(m+p),$ then
\[D_{\Phi,k}(d^k;\varrho_k)\leq-\gamma_\Phi\|d^k\|^2.\]
\end{proposition}
\begin{proof}
Taking the inner product of the approximate QCQP stationarity equation with $d^k$ and substituting the quadratic constraint models gives
\[\begin{aligned} (g^k)^\top d^k
={}&-(d^k)^\top G_\zeta^kd^k+\sum_i\widehat\lambda_i^k
\bigl(h_i(x^k)-\widehat h_i^k(d^k)\bigr)+\sum_j\widehat\mu_j^k
\bigl(c_j(x^k)-\widehat c_j^k(d^k)\bigr)\\
&-\frac12\sum_i\widehat\lambda_i^k
(d^k)^\top C_{h_i}^kd^k-\frac12\sum_j\widehat\mu_j^k
(d^k)^\top C_{c_j}^kd^k+(u^k)^\top d^k.
\end{aligned}\]

Using the model residual bounds, multiplier bounds, and symmetry of the curvature matrices yields
\[\begin{aligned}(g^k)^\top d^k\leq{}&-(\zeta_k-M_G-\tfrac12\ell_k M_C^\Sigma)\|d^k\|^2\\&+\ell_k V(x^k)
+\delta_k(\|d^k\|+m\ell_k+p).\end{aligned}\]

The one-sided directional derivative of the exact $L_1$ constraint violation satisfies
\[\begin{aligned} V'(x^k;d^k)\leq{}& -V(x^k)+(m+p)\delta_k+\frac12M_C^\Sigma\|d^k\|^2.\end{aligned}\]

Combining these inequalities proves the first claim.
Finally, from \eqref{eq:relative_qcqp_accuracy},
\[\begin{aligned}\delta_k(\|d^k\|+A_k)&\leq\kappa_\delta
\frac{\|d^k\|^2}{1+\|d^k\|}(\|d^k\|+A_k)\\
&\leq\kappa_\delta\max\{1,A_k\}\|d^k\|^2.\end{aligned}\]

The stated regularization condition and $\varrho_k\geq\ell_k$ imply the conclusion.
\end{proof}

\subsection{Merit Summability and PAGE Error Control}
\begin{lemma}[Summability of accepted search directions]
\label{lem:merit_summability}
Suppose Assumption~\ref{ass:global} holds. Then
\begin{equation}\label{eq:descent_summability}
\sum_{k=0}^{\infty}\alpha_k\|d^k\|^2<\infty, \end{equation}
and, with $s^k=\alpha_kd^k$,
\[\sum_{k=0}^{\infty}\|s^k\|^2<\infty.\]
Moreover, both sums admit finite deterministic upper bounds under the constants in Assumption~\ref{ass:global}.
\end{lemma}
\begin{proof}
By the accepted sufficient-descent and Armijo conditions,
\[ \Phi(x^{k+1};\varrho_k)\leq\Phi(x^k;\varrho_k)
-\sigma\gamma_\Phi \alpha_k\|d^k\|^2. \]

For consecutive accepted trial penalties,
\[ \begin{aligned} \Phi(x^{k+1};\varrho_{k+1})={}&
\Phi(x^{k+1};\varrho_k)+(\varrho_{k+1}-\varrho_k) V(x^{k+1}).
\end{aligned} \]

Compactness and continuity imply the existence of $M_V<\infty$ such that $V(x^k)\leq M_V$.

Consequently,
\[ \begin{aligned} \Phi(x^{k+1};\varrho_{k+1}) \leq{}&
\Phi(x^k;\varrho_k) -\sigma\gamma_\Phi \alpha_k\|d^k\|^2+M_V|\varrho_{k+1}-\varrho_k|. \end{aligned}\]

Summing from $k=0$ to $K-1$ gives
\[ \begin{aligned}
\sigma\gamma_\Phi\sum_{k=0}^{K-1}\alpha_k\|d^k\|^2
\leq{}&\Phi(x^0;\varrho_0) -\Phi(x^K;\varrho_K)+M_V\sum_{k=0}^{K-1} |\varrho_{k+1}-\varrho_k|.\end{aligned} \]

Since $\varrho_K>0$ and $V(x^K)\geq0$,
\[ \Phi(x^K;\varrho_K)\geq f_{\min}. \]

Using the bounded trial-penalty variation yields
\[ \sum_{k=0}^{K-1} \alpha_k\|d^k\|^2 \leq \frac{ \Phi(x^0;\varrho_0)-f_{\min} +M_VB_{\mathrm{tr}} }{
\sigma\gamma_\Phi }. \]

Because $\varrho_0\leq\overline\varrho$, the right-hand side can be bounded by a deterministic constant $T_d$. Letting $K\to\infty$ proves \eqref{eq:descent_summability}.

Finally, $0<\alpha_k\leq1$ implies
\[ \|s^k\|^2 =\alpha_k^2\|d^k\|^2 \leq\alpha_k\|d^k\|^2, \]
which proves the second assertion.
\end{proof}

\begin{lemma}[Summability and almost-sure decay of PAGE errors]
\label{lem:page_error}
Suppose Assumption~\ref{ass:global} holds, and $\mathbb E\|e^0\|^2<\infty$.

Then
\[ \sum_{k=0}^{\infty} \mathbb E\|e^k\|^2 \leq T_e<\infty, \]
where one may take
$ T_e=\frac{ \mathbb E\|e^0\|^2+ \frac{(1-p_{\mathrm{ref}})L_F^2}{b_s}T_d }{p_{\mathrm{ref}}}. $

Consequently, \[ e^k\longrightarrow0 \qquad\text{almost surely}. \]
\end{lemma}
\begin{proof}
Taking expectations in \eqref{eq:page_error_recursion} gives
\[ \begin{aligned}
\mathbb E\|e^{k+1}\|^2 \leq{}& (1-p_{\mathrm{ref}})
\mathbb E\|e^k\|^2+\frac{(1-p_{\mathrm{ref}})L_F^2}{b_s}
\mathbb E\|s^k\|^2. \end{aligned} \]

Summing this recursion and using the deterministic bound from Lemma~\ref{lem:merit_summability} establishes
\[ p_{\mathrm{ref}} \sum_{k=0}^{K}\mathbb E\|e^k\|^2
\leq \mathbb E\|e^0\|^2+ \frac{(1-p_{\mathrm{ref}})L_F^2}{b_s}T_d. \]

The upper bound is independent of $K$. Thus,
\[ \sum_{k=0}^{\infty} \mathbb E\|e^k\|^2<\infty. \]

By Tonelli's theorem, $\sum_k\|e^k\|^2<\infty$ almost surely.
Hence $\|e^k\|\to0$ almost surely.
\end{proof}

\subsection{Transfer of QCQP Optimality to NLP KKT Conditions}
Recall the KKT residual mapping $\mathcal R(x,\lambda,\mu)$ and its Euclidean norm $\mathcal E(x,\lambda,\mu)$ defined
in~\eqref{eq:kkt_residual}. For the convergence analysis, define $\mathcal E_k=\mathcal E\left(x^k,\widehat\lambda^k,
\widehat\mu^k\right),$ where $(\widehat\lambda^k,\widehat\mu^k)$ are the multipliers recovered from the accepted QCQP model.

\begin{lemma}[QCQP-to-NLP KKT transfer]\label{lem:kkt_transfer}
Suppose (G1)--(G2) hold and \eqref{eq:qcqp_accuracy} is satisfied. Then there exists $K_R>0$, independent of $k$,
such that
\begin{equation} \label{eq:kkt_transfer_bound}
\mathcal   E_k\leq K_R\left(\|d^k\|+\|d^k\|^2+\|e^k\|+\delta_k
\right).\end{equation}\end{lemma}
\begin{proof}
The approximate QCQP stationarity equation implies
\[ \begin{aligned}
\nabla_x\mathcal L(x^k,\widehat\lambda^k,\widehat\mu^k)
={}&u^k-e^k-G_\zeta^kd^k-\sum_i\widehat\lambda_i^k
C_{h_i}^kd^k-\sum_j\widehat\mu_j^k C_{c_j}^kd^k. \end{aligned} \]

Uniform curvature and multiplier bounds therefore give
\[ \|\nabla_x\mathcal L(x^k,\widehat\lambda^k,\widehat\mu^k)\|
\leq K_1\|d^k\|+\|e^k\|+\delta_k. \]

For each equality constraint,
\[ h_i(x^k) = \widehat h_i^k(d^k) -\nabla h_i(x^k)^\top d^k
-\frac12(d^k)^\top C_{h_i}^kd^k. \]

Hence \[ \|h(x^k)\| \leq K_2\left( \|d^k\|+\|d^k\|^2+\delta_k
\right). \]

The same argument applied to the positive parts of the inequality constraints yields \[ \|[c(x^k)]^+\| \leq K_3\left( \|d^k\|+\|d^k\|^2+\delta_k \right). \]

Furthermore, bounded multipliers and QCQP complementarity imply
\[ \|\widehat\mu^k\odot c(x^k)\| \leq K_4\left( \|d^k\|+\|d^k\|^2+\delta_k \right). \]

Finally, the dual-feasibility residual is bounded by $\delta_k$; it vanishes when the recovered inequality multipliers are explicitly nonnegative. Combining these estimates proves the assertion.
\end{proof}

\begin{theorem}[Almost-sure subsequential KKT convergence]
\label{thm:subsequential_kkt}
Suppose Assumptions~\ref{ass:global} and \ref{ass:relative_qcqp} hold.

If $\sum_{k=0}^{\infty}\alpha_k=\infty, $ then
\[ \liminf_{k\to\infty}\|d^k\|=0 \qquad\text{almost surely},\]
and the generated primal sequence possesses at least one KKT accumulation point almost surely.
\end{theorem}

\begin{proof}
Lemma~\ref{lem:merit_summability} gives $\sum_{k=0}^{\infty}\alpha_k\|d^k\|^2<\infty.$ If $\liminf_k\|d^k\|>0$, then $\|d^k\|\geq\varepsilon$ for all sufficiently large $k$, for some $\varepsilon>0$. This would imply
\[ \sum_k\alpha_k\|d^k\|^2 \geq\varepsilon^2 \sum_k\alpha_k=\infty, \] a contradiction.

Thus there exists a subsequence $\{k_\ell\}$ such that
$d^{k_\ell}\to0.$ By Lemma~\ref{lem:page_error}, $e^{k_\ell}\to0$ almost surely. The relative QCQP accuracy condition implies $\delta_{k_\ell}\to0.$ It follows from Lemma~\ref{lem:kkt_transfer} that $ \mathcal R_{k_\ell}\to0 \qquad\text{almost surely}. $

Compactness of the primal iterates and boundedness of the local multipliers permit extraction of a further subsequence for which
\[ (x^{k_\ell}, \widehat\lambda^{k_\ell}, \widehat\mu^{k_\ell})
\to(x^\star,\lambda^\star,\mu^\star). \]

By continuity of the KKT residual components, the limit satisfies \eqref{eq:kkt} Therefore $x^\star$ is a KKT point.
\end{proof}

\begin{corollary}[Full KKT residual convergence] \label{cor:full_kkt}
Under the assumptions of Theorem~\ref{thm:subsequential_kkt},
suppose there exist $k_0$ and $\underline\alpha>0$ such that
$ \alpha_k\geq\underline\alpha, \  k\geq k_0. $

Then
\[\|d^k\|\to0,\qquad \delta_k\to0,\qquad \mathcal R_k\to0
\ \ \ almost \  surely. \]

Every primal accumulation point is a KKT point.
\end{corollary}

\begin{proof}
The uniform step-length bound and
\eqref{eq:descent_summability} imply
\[ \sum_{k=k_0}^{\infty}\|d^k\|^2\leq\frac{T_d}{\underline\alpha} <\infty. \]

Hence $\|d^k\|\to0$. Assumption~\ref{ass:relative_qcqp}
implies $\delta_k\to0$, while Lemma~\ref{lem:page_error} gives
$e^k\to0$ almost surely. The conclusion now follows from
\eqref{eq:kkt_transfer_bound} and continuity of the nonlinear-programming KKT conditions.
\end{proof}

\subsection{Step-Length Conditions and Complexity}

For a nonzero direction $d^k$, let $\tau_k\in(0,1]$ denote the first positive zero crossing, within $[0,1]$, of any affine
constraint prediction $r(x^k)+\alpha\nabla r(x^k)^\top d^k,$
with $\tau_k=1$ when no such crossing occurs.

\begin{lemma}[Finite Armijo acceptance] \label{lem:finite_armijo_inexact}
Suppose the accepted candidate direction satisfies
$D_{\Phi,k}(d^k;\varrho_k) \leq-\gamma_\Phi\|d^k\|^2.$
Assume the relevant gradients are Lipschitz
continuous, $\|e^k\|\leq\theta\gamma_\Phi\|d^k\|,$ and $0\leq\theta<1-\sigma.$

Then the exact-merit Armijo condition is satisfied for all sufficiently small positive trial steps. In exact arithmetic, backtracking terminates after finitely many reductions.

More precisely, for a suitable $L_{\Phi,k}$, backtracking accepts a step satisfying
\[ \alpha_k\geq \beta\min\left\{ 1,\tau_k, \frac{ 2(1-\sigma-\theta)\gamma_\Phi }{ L_{\Phi,k} } \right\}, \]
where a zero value of $L_{\Phi,k}$ is interpreted through an infinite third bound.
\end{lemma}
\begin{proof}
On $0\leq\alpha\leq\tau_k$, the absolute-value and positive-part functions applied to the affine constraint predictions remain on their initial linear branches.

Lipschitz continuity of the gradients therefore gives
\[ \begin{aligned}
\Phi(x^k+\alpha d^k;\varrho_k) \leq{}&
\Phi(x^k;\varrho_k)+\alpha D_{\Phi,k}(d^k;\varrho_k)\\
&+\alpha\|e^k\|\|d^k\|+\frac12L_{\Phi,k} \alpha^2\|d^k\|^2.
\end{aligned} \]

By sufficient descent and the stated PAGE relative-error bound,
\[ \begin{aligned}
&\Phi(x^k+\alpha d^k;\varrho_k) -\Phi(x^k;\varrho_k)-\sigma\alpha D_{\Phi,k}(d^k;\varrho_k)\\ &\leq
-\alpha(1-\sigma-\theta) \gamma_\Phi\|d^k\|^2
+\frac12L_{\Phi,k} \alpha^2\|d^k\|^2. \end{aligned} \]

The right-hand side is nonpositive whenever
\[ \alpha\leq \frac{ 2(1-\sigma-\theta)\gamma_\Phi }{
L_{\Phi,k} }. \]

Combining this requirement with $\alpha\leq\tau_k$ proves finite acceptance and the backtracking bound.
\end{proof}
\begin{remark}
The preceding lemma ensures finite Armijo acceptance but does not by itself imply a uniform positive lower bound
on the accepted step lengths. If, in addition, $ \tau_k\geq\underline\tau>0$ and $ L_{\Phi,k}\leq\overline L_\Phi<\infty $ for all sufficiently large $k$, and the
relative PAGE error condition holds with fixed $\theta<1-\sigma$, then
\[ \alpha_k\geq \beta\min\left\{ 1,\underline\tau,
\frac{ 2(1-\sigma-\theta)\gamma_\Phi }{ \overline L_\Phi
} \right\} =:\underline\alpha>0. \]
Thus, the uniform step-length hypothesis used in the subsequent convergence and complexity results follows under these additional conditions.
\end{remark}

\begin{theorem}[Expected KKT complexity] \label{thm:expected_kkt_complexity}
Suppose Assumptions~\ref{ass:global} and
\ref{ass:relative_qcqp} hold from some deterministic index $k_0$ onward. Assume that $\alpha_k\geq\underline\alpha>0,
\ \  k\geq k_0,$ where $\underline\alpha$ is deterministic.
Then there exists a deterministic constant
$C_{\mathrm{rate}}<\infty$ such that
\[
\min_{k_0\leq k<K} \mathbb E[\mathcal E^2] \leq
\frac{C_{\mathrm{rate}}}{K-k_0}.\]
If $R$ is chosen uniformly from $\{k_0,\ldots,K-1\}$, independently of the algorithmic randomness, then
\[ \mathbb E[\mathcal E_R^2] \leq \frac{C_{\mathrm{rate}}}{K-k_0}. \]
Consequently, an expected squared KKT residual of order $\varepsilon^2$ requires $O(\varepsilon^{-2})$ outer iterations
under these assumptions.
\end{theorem}

\begin{proof}
By Lemma~\ref{lem:merit_summability},
$\sum_{k=k_0}^{\infty}\|d^k\|^2\leq\frac{T_d}{\underline\alpha}.$ In particular, the directions satisfy the
deterministic bound
$ \|d^k\|\leq M_d:= \sqrt{T_d/\underline\alpha},
\ \  k\geq k_0. $ Assumption~\ref{ass:relative_qcqp} gives
$\delta_k\leq\kappa_\delta\|d^k\|.$ Applying Lemma~\ref{lem:kkt_transfer} and using the direction bound yields
\[\mathcal E_k\leq K_R\left[ (1+M_d+\kappa_\delta)\|d^k\|
+\|e^k\|\right].\]

Therefore,
\[ \begin{aligned}
\mathbb E[\mathcal E_k^2] \leq{}&
2K_R^2(1+M_d+\kappa_\delta)^2
\mathbb E\|d^k\|^2+2K_R^2\mathbb E\|e^k\|^2.
\end{aligned} \]

Summing over $k=k_0,\ldots,K-1$, and using Lemmas~\ref{lem:merit_summability} and \ref{lem:page_error}, gives
\[ \sum_{k=k_0}^{K-1} \mathbb E[\mathcal E_k^2]
\leq C_{\mathrm{rate}}, \]
where, for example,
\[ \begin{aligned}
C_{\mathrm{rate}} =2K_R^2\Bigg[ &(1+M_d+\kappa_\delta)^2
\frac{T_d}{\underline\alpha} +T_e \Bigg]. \end{aligned} \]
Dividing by $K-k_0$ yields the bound on the average expected squared KKT residual. The minimum is no larger than this average. Uniform independent selection of $R$ gives
the same average directly, completing the proof.
\end{proof}
\begin{remark}[Complexity interpretation]
The complexity result counts accepted outer iterations, not component-gradient evaluations, QCQP inner iterations, regularization re-solves, or line-search function evaluations.
Consequently, it should not be interpreted as an $O(\varepsilon^{-2})$ stochastic-oracle complexity
guarantee. The constants may depend on the trial-penalty
variation bound, sampling parameters, QCQP accuracy, and the assumed lower bound on accepted step lengths.
\end{remark}

\subsection{Local Consistency of Learned Constraint Curvature}

The preceding global results do not require the constraint SR1 matrices to converge to the exact constraint Hessians in operator norm. To examine local second-order modeling, consider
a twice continuously differentiable constraint function $r$, with Lipschitz-continuous Hessian. Suppose a curvature update satisfies the secant equation
\[ C_r^{k+1}s^k=\nabla r(x^{k+1})-\nabla r(x^k),\qquad
s^k=x^{k+1}-x^k.\]

\begin{proposition}[Secant-direction consistency]
\label{prop:sr1_secant_consistency}
If $\nabla^2r$ is Lipschitz continuous with constant $L_{H,r}$, then every update satisfying the preceding secant equation obeys
\[ \left\| \left( C_r^{k+1}-\nabla^2r(x^{k+1}) \right)s^k
\right\| \leq \frac{L_{H,r}}2\|s^k\|^2. \]
\end{proposition}
\begin{proof}
The fundamental theorem of calculus gives
\[ \nabla r(x^{k+1})-\nabla r(x^k) = \int_0^1
\nabla^2r(x^k+ts^k)s^k\,dt. \]

Subtracting $\nabla^2r(x^{k+1})s^k$ and using Hessian Lipschitz continuity yields
\[ \begin{aligned}
&\left\| \left( C_r^{k+1}-\nabla^2r(x^{k+1}) \right)s^k
\right\|\\ &\leq \int_0^1 L_{H,r}(1-t) \|s^k\|^2\,dt\\
&= \frac{L_{H,r}}2\|s^k\|^2. \end{aligned} \] \end{proof}
\begin{assumption}[Directional curvature consistency]
\label{ass:direction_consistency}
Along a sequence approaching a KKT point,
the learned constraint models satisfy
\[
\left\|
(C_r^k-\nabla^2r(x^k))d^k
\right\|
=o(\|d^k\|)
\]
for each nonlinear constraint $r$.
\end{assumption}

\begin{lemma}[Second-order constraint-model consistency]
\label{lem:constraint_model_consistency}
Under Assumption~\ref{ass:direction_consistency}
and Lipschitz continuity of the constraint Hessians, the learned quadratic model satisfies \[ r(x^k+d^k)-\widehat r^k(d^k) =o(\|d^k\|^2) \ \ as \ \  \|d^k\|\to0.\] \end{lemma}
\begin{proof}
Taylor's theorem gives
\[ \begin{aligned}
r(x^k+d^k) ={}&r(x^k)+\nabla r(x^k)^\top d^k+\frac12(d^k)^\top
\nabla^2r(x^k)d^k +O(\|d^k\|^3). \end{aligned} \]

Subtracting the learned quadratic model yields
\[ \begin{aligned}
r(x^k+d^k)-\widehat r^k(d^k) ={}& \frac12(d^k)^\top \left(
\nabla^2r(x^k)-C_r^k \right)d^k+O(\|d^k\|^3). \end{aligned} \]

By directional consistency, \[ \left| (d^k)^\top (\nabla^2r(x^k)-C_r^k)d^k \right| =o(\|d^k\|^2). \]

The result follows.
\end{proof}

\begin{theorem}[Local stability relative to exact SQCQP]
\label{thm:direction_stability}
Consider, at the same base point $x^k$, the learned-curvature QCQP and the QCQP using the exact constraint Hessians. Suppose both problems use the same stochastic gradient, objective curvature model, and objective regularization.
Let $d_{\mathrm{QN}}^k$ and $d_{\mathrm{EX}}^k$ denote local KKT directions on corresponding nearby active-set branches.
Assume the relevant branches are locally well defined, their multipliers remain bounded, and the KKT systems satisfy a uniform local inverse error bound. Such a bound may be obtained under suitable strong regularity, active-set stability, and second-order sufficient conditions. Suppose additionally that
\[ \left\| (C_r^k-\nabla^2r(x^k)) d_{\mathrm{EX}}^k \right\|
=o(\|d_{\mathrm{EX}}^k\|) \]
for every nonlinear constraint $r$. Then, as the nonzero exact comparison directions approach zero,
\begin{equation} \label{eq:direction_stability}
\|d_{\mathrm{QN}}^k-d_{\mathrm{EX}}^k\| =o(\|d_{\mathrm{EX}}^k\|). \end{equation} \end{theorem}
\begin{proof}
Evaluate the learned-curvature QCQP KKT equations at the exact-Hessian QCQP KKT pair. Because both subproblems use the same objective gradient and objective curvature model, the
stationarity residual difference arises from the constraint curvature terms:
$ \sum_i\widehat\lambda_i^k (C_{h_i}^k-\nabla^2h_i(x^k))
d_{\mathrm{EX}}^k $ together with the corresponding inequality terms. Bounded multipliers and directional curvature
consistency imply that the resulting stationarity perturbation is $o(\|d_{\mathrm{EX}}^k\|)$. The differences between the quadratic constraint models at the exact-Hessian direction are
\[ \frac12(d_{\mathrm{EX}}^k)^\top (C_r^k-\nabla^2r(x^k))
d_{\mathrm{EX}}^k = o(\|d_{\mathrm{EX}}^k\|^2). \]
The associated active constraint and complementarity perturbations therefore have higher order than $\|d_{\mathrm{EX}}^k\|$. The assumed uniform local inverse error bound for the corresponding KKT branches converts
these KKT residual perturbations into a primal-direction perturbation of order $o(\|d_{\mathrm{EX}}^k\|)$.
This establishes \eqref{eq:direction_stability}.
\end{proof}

\section{Numerical Experiments}
\label{sec:experiments}

The numerical performance of the proposed VR-QN-SQCQP method is evaluated on two constrained finite-sum optimization problems: sparse principal component analysis (Sparse PCA) using the CDC BRFSS 2015 dataset and a dense nonlinear constrained benchmark. A comparison is performed with
Exact-SQCQP-PAGE, in which exact constraint Hessians
are employed instead of the learned SR1 constraint-curvature approximations. Identical PAGE gradient estimators, objective SR1 updates, QCQP solution procedures, and globalization
mechanisms are employed in both methods.

The experiments are conducted to assess convergence reliability, KKT accuracy, computational runtime, gradient-oracle complexity, QCQP workload, and constraint-curvature construction cost. Both methods are implemented in Python using NumPy and SciPy, and the QCQP subproblems are solved
using SLSQP.

For both experiments, the PAGE small-batch size and refresh probability are fixed at $b_{\mathrm{s}}=256$ and $p_{\mathrm{ref}}=0.40$, respectively. Full-gradient evaluations are employed during PAGE refreshes. The objective
SR1 safeguard is selected as $\kappa_G=5\times10^{-3}$, while the constraint SR1 denominator safeguard is fixed at
$\kappa_C=10^{-8}$. The Armijo parameters are set to $\sigma=10^{-4}$ and $\beta=0.5$. A maximum of 300 outer iterations is permitted, and convergence is declared when the KKT residual satisfies $\mathcal{E}(x,\lambda,\mu)\leq10^{-5}$.
Ten random seeds are considered for each benchmark, with identical problem instances and shared algorithmic settings employed for the compared methods.

The first experiment is performed on the full CDC BRFSS 2015 dataset containing $N=253\,680$ observations. A total of $n=20$
predictors is selected using a fixed feature-selection seed. The Sparse PCA problem is formulated as
$\begin{aligned}
\min_{x\in\mathbb{R}^{n}}  \ -x^\top\Sigma x \ \ \ \
\mathrm{s.\ t.} \ \ \|x\|^2-1=0, \ \ \
\sum_{i=1}^{n}
\left(
\sqrt{x_i^2+\varepsilon^2}
-\varepsilon
\right)-\tau\leq0,
\end{aligned}$
where $\Sigma$ denotes the sample covariance matrix, $\tau=3$, and $\varepsilon=10^{-3}$. The initial point is selected as $x^0=e_1$. For this problem, the exact equality-constraint
Hessian is $2I$, while the inequality-constraint Hessian is diagonal, with diagonal entries $\frac{\varepsilon^2}{(x_i^2+\varepsilon^2)^{3/2}}.$ The constraint-SR1 Frobenius-norm safeguard is selected as $M_c=4472.136$, based on
the analytical bound $\|\nabla^2c(x)\|_F\leq\sqrt{n}/\varepsilon$.

The numerical comparisons are presented in Table~\ref{tab:final_results} and Fig.~\ref{fig:sparse_results}. Convergence to the prescribed KKT tolerance is achieved by both methods on all ten seeds.
The median numbers of outer iterations are 66 and 41 for VR-QN-SQCQP and Exact-SQCQP-PAGE, respectively.

Essentially identical final objective values of approximately $-3.11025908$ are obtained by both methods. The explained variance is approximately $3.110259$, corresponding to
a variance-retention ratio of $90.496\%$, while a support size of 15 is reported. Negligible constraint violations are observed, and the final KKT residuals remain below $10^{-5}$ for all seeds. The median computational runtimes are $3.384$~s and $1.747$~s for the proposed and exact-Hessian methods, respectively. The corresponding median QCQP solve counts
are 400 and 203, while the median gradient-oracle sample counts are 9\,658\,272 and 4\,960\,584. Greater computational efficiency is therefore obtained by Exact-SQCQP-PAGE
on this structured benchmark. Nevertheless, comparable solution quality is achieved by VR-QN-SQCQP without any exact constraint-Hessian evaluations.

The effect of the constraint-SR1 safeguard is also illustrated in Fig.~\ref{fig:sparse_results}. With the original Frobenius-norm cap of 1000, none of the ten proposed-method runs reaches the KKT tolerance within 1000 iterations. After the curvature-informed cap of 4472.136 is employed, convergence is obtained on all ten seeds within 131 iterations. The importance of
selecting curvature safeguards consistent with the nonlinear-constraint scale is thereby demonstrated. The two safeguard configurations are compared using their respective iteration budgets.

The second experiment is conducted on a dense nonlinear constrained optimization problem with $N=20\,000$ samples,
$n=100$ variables, and $q=2000$ nonlinear constraint terms.

The finite-sum objective is defined as
$f(x)=\frac{1}{2N}\|A_{\mathrm{obj}}x-y\|^2+\frac{\gamma}{2}\|x\|^2,\qquad\gamma=10^{-3}.$ One nonlinear equality constraint and one nonlinear inequality constraint are
constructed using $\phi_A(x)=\frac{1}{q}\sum_{\ell=1}^{q}
\log\cosh(a_\ell^\top x).$ The corresponding constraints are
$h(x)=\phi_{A_h}(x)-b_h=0$ and $c(x)=\phi_{A_c}(x)-b_c\leq0$,
where the offsets are generated using a feasible reference point. The instance generation seed is fixed at 2026. The constraint-SR1 Frobenius-norm cap is selected as 100.
The exact Hessian of the nonlinear constraint-building function is $\nabla^2\phi_A(x)=\frac{1}{q}A^\top\operatorname{diag}
\left(\operatorname{sech}^2(Ax)\right)A.$

These constraint Hessians are generally dense and state dependent. Their explicit construction requires $O(qn^2)$ operations, whereas an SR1 curvature update requires
$O(n^2)$ operations after constraint-gradient evaluation. A computational benefit from learned constraint curvature is therefore expected when the nonlinear constraint dimension is large.

The results are summarized in Table~\ref{tab:final_results} and
Fig.~\ref{fig:dense_results}. Convergence to the prescribed tolerance is achieved by VR-QN-SQCQP on all ten seeds, while Exact-SQCQP-PAGE converges on nine. Median computational statistics are calculated over the nine jointly converged seeds.

The median numbers of outer iterations are 32 and 40, respectively, while the corresponding median QCQP solve
counts are 83 and 130. The median gradient-oracle sample counts are 389\,728 and 491\,264. Median runtimes of approximately $2.695$~s and $2.920$~s are recorded. A more pronounced difference is observed in constraint-curvature construction time. The median cumulative times are approximately $0.0111$~s for VR-QN-SQCQP and $0.1230$~s for
Exact-SQCQP-PAGE, indicating an approximately elevenfold reduction in this computational component. For seed 9, the proposed method reaches the prescribed KKT tolerance in 64 iterations, whereas the exact-Hessian method reaches the 300-iteration limit with a final KKT residual of approximately $5.74\times10^{-5}$. The final objective values obtained
by the two methods remain close to $0.107014$. Across the ten dense-benchmark runs, zero exact constraint-Hessian
evaluations are recorded for the proposed method, compared with 1350 evaluations for Exact-SQCQP-PAGE.

The results from both experiments demonstrate that safeguarded SR1 curvature models can be incorporated into stochastic SQCQP
without explicit constraint-Hessian evaluations while retaining comparable solution quality. For inexpensive structured constraint Hessians, the exact-Hessian approach remains computationally effective. For dense nonlinear constraints, substantial reductions in curvature-construction cost are obtained through the proposed method.

\begin{table*}[h]
\centering
\caption{Numerical comparison of VR-QN-SQCQP
and Exact-SQCQP-PAGE. Sparse PCA medians are
computed over ten converged seeds, while dense
nonlinear medians are computed over nine jointly
converged seeds. Convergence counts and total
constraint-Hessian evaluations include all
ten seeds.}
\label{tab:final_results}
\small

\begin{tabular}{lcccc}
\hline
&
\multicolumn{2}{c}{\textbf{Sparse PCA}}
&
\multicolumn{2}{c}{\textbf{Dense nonlinear}}\\

\textbf{Performance measure}
& \textbf{VR-QN}
& \textbf{Exact}
& \textbf{VR-QN}
& \textbf{Exact}
\\
\hline

Converged runs
& 10/10 & 10/10 & 10/10 & 9/10
\\

Median outer iterations
& 66 & 41 & 32 & 40
\\

Median runtime (s)
& 3.384 & 1.747 & 2.695 & 2.920
\\

Median QCQP solves
& 400 & 203 & 83 & 130
\\

Median gradient samples
& 9\,658\,272
& 4\,960\,584
& 389\,728
& 491\,264
\\

Median curvature time (s)
& -- & -- & 0.0111 & 0.1230
\\

Exact constraint Hessians (total)
& 0 & 866 & 0 & 1350
\\

\hline
\end{tabular}
\end{table*}

\newcommand{\paperfigpanel}[2]{%
\begin{minipage}[t]{0.485\textwidth}
\centering
\includegraphics[
    width=\linewidth,
    height=0.165\textheight,
    keepaspectratio
]{#1}
\par\vspace{1pt}
{\scriptsize #2}
\end{minipage}%
}

\clearpage

\begin{figure*}[p]
\centering

\paperfigpanel{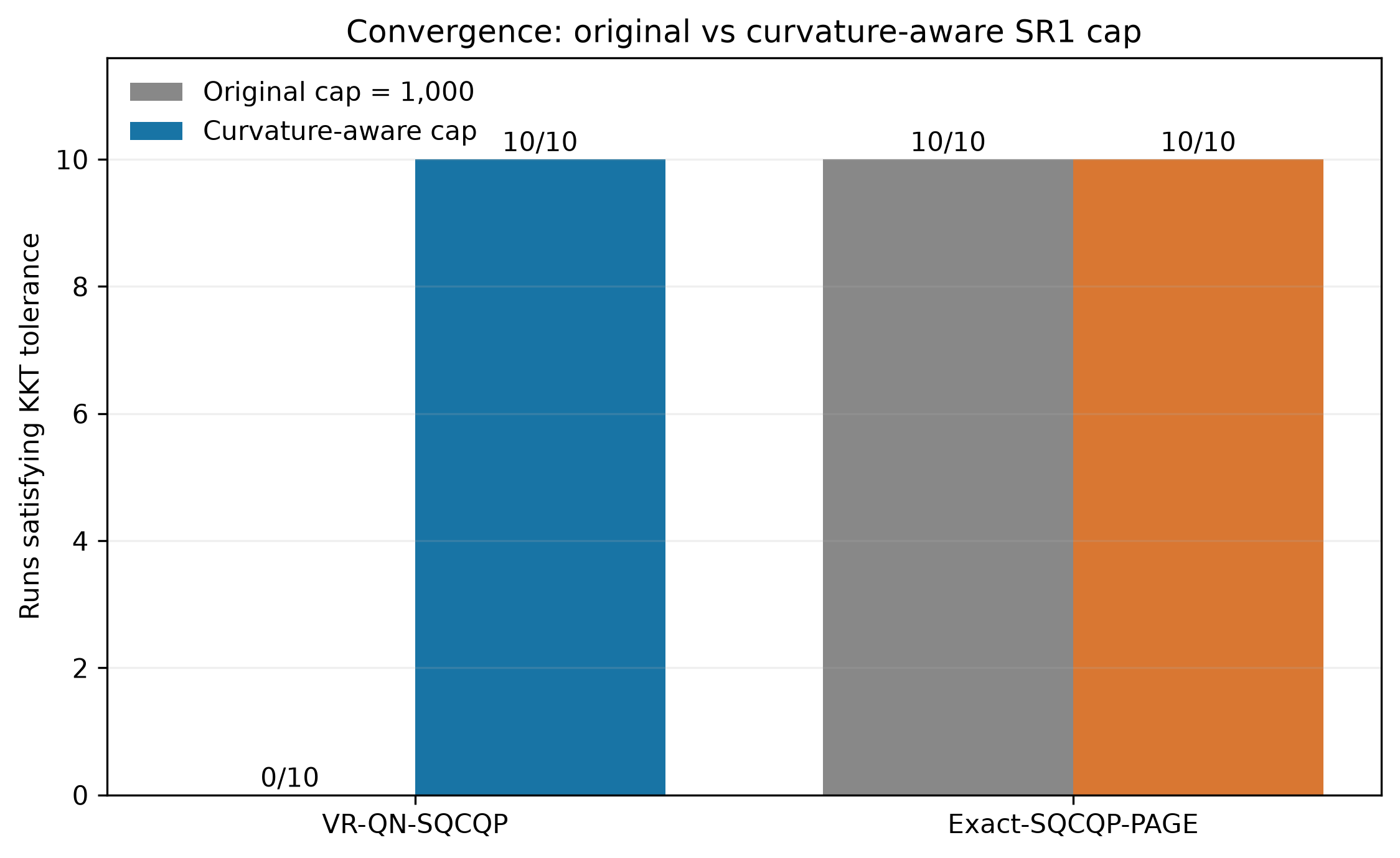}{\textbf{(a)} Convergence reliability.}\hfill
\paperfigpanel{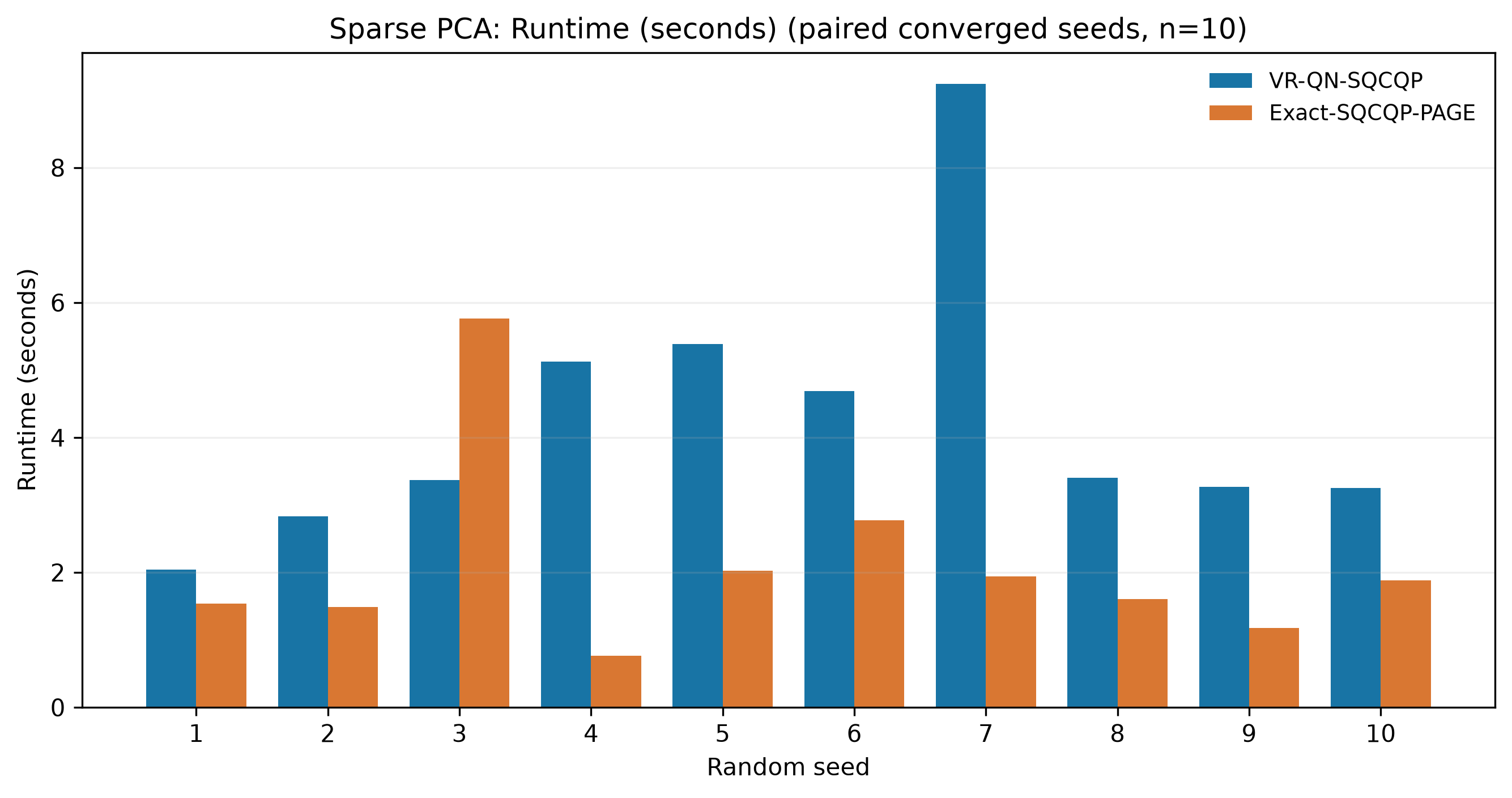}{\textbf{(b)} Computational runtime.}

\par\vspace{3mm}

\paperfigpanel{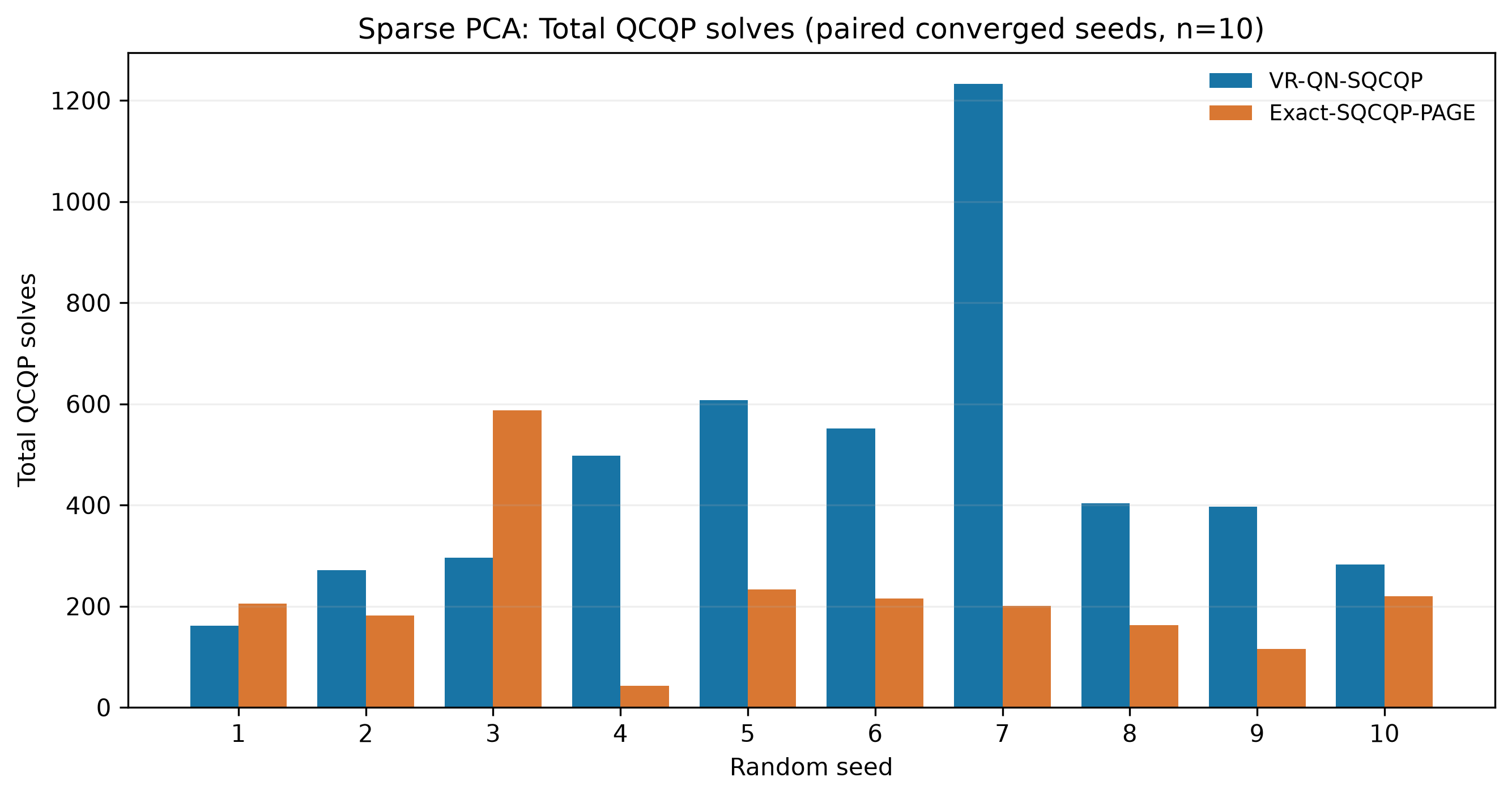}{\textbf{(c)} Total QCQP solves.}\hfill
\paperfigpanel{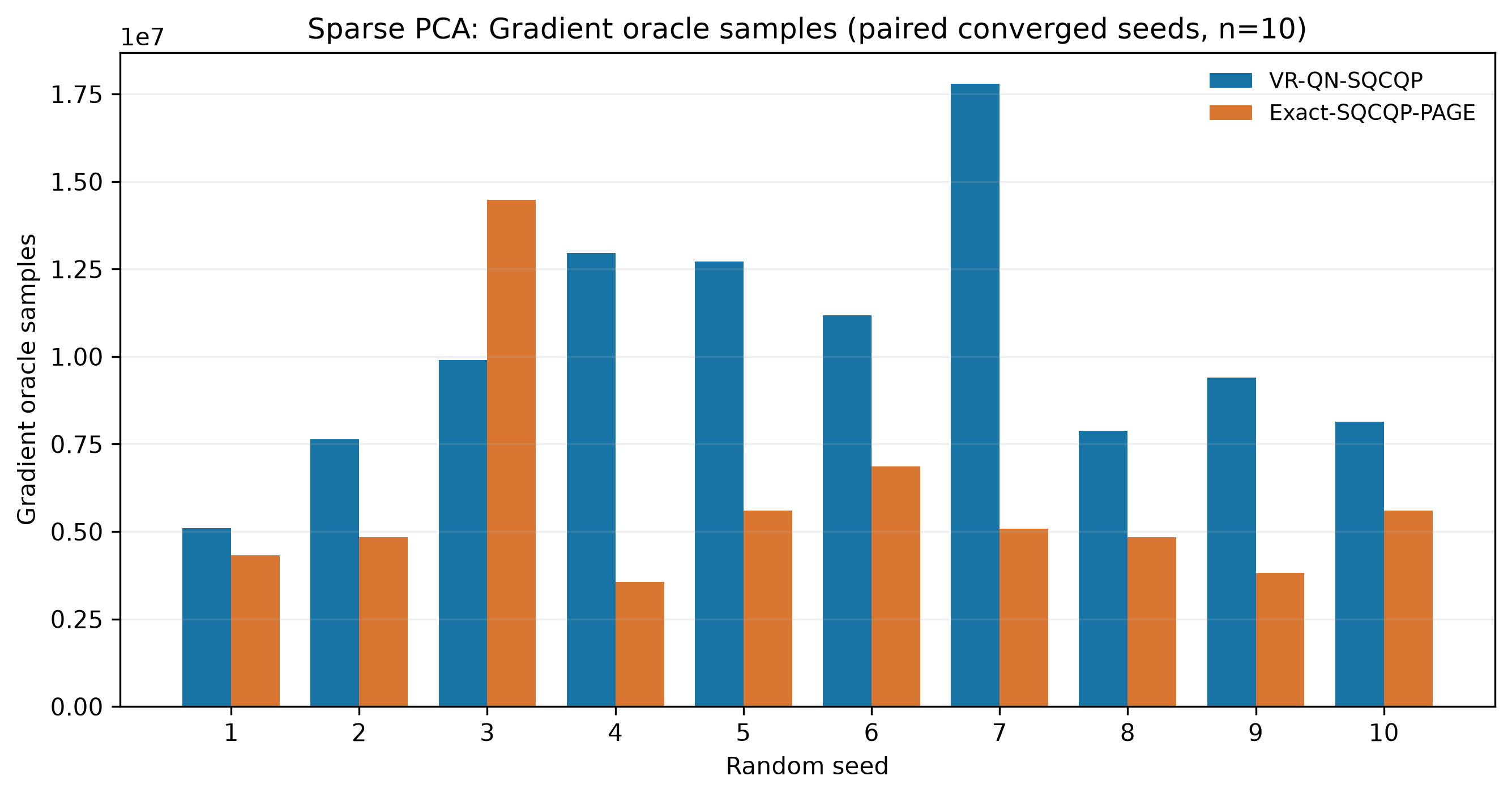}{\textbf{(d)} Gradient-oracle samples.}

\par\vspace{3mm}

\paperfigpanel{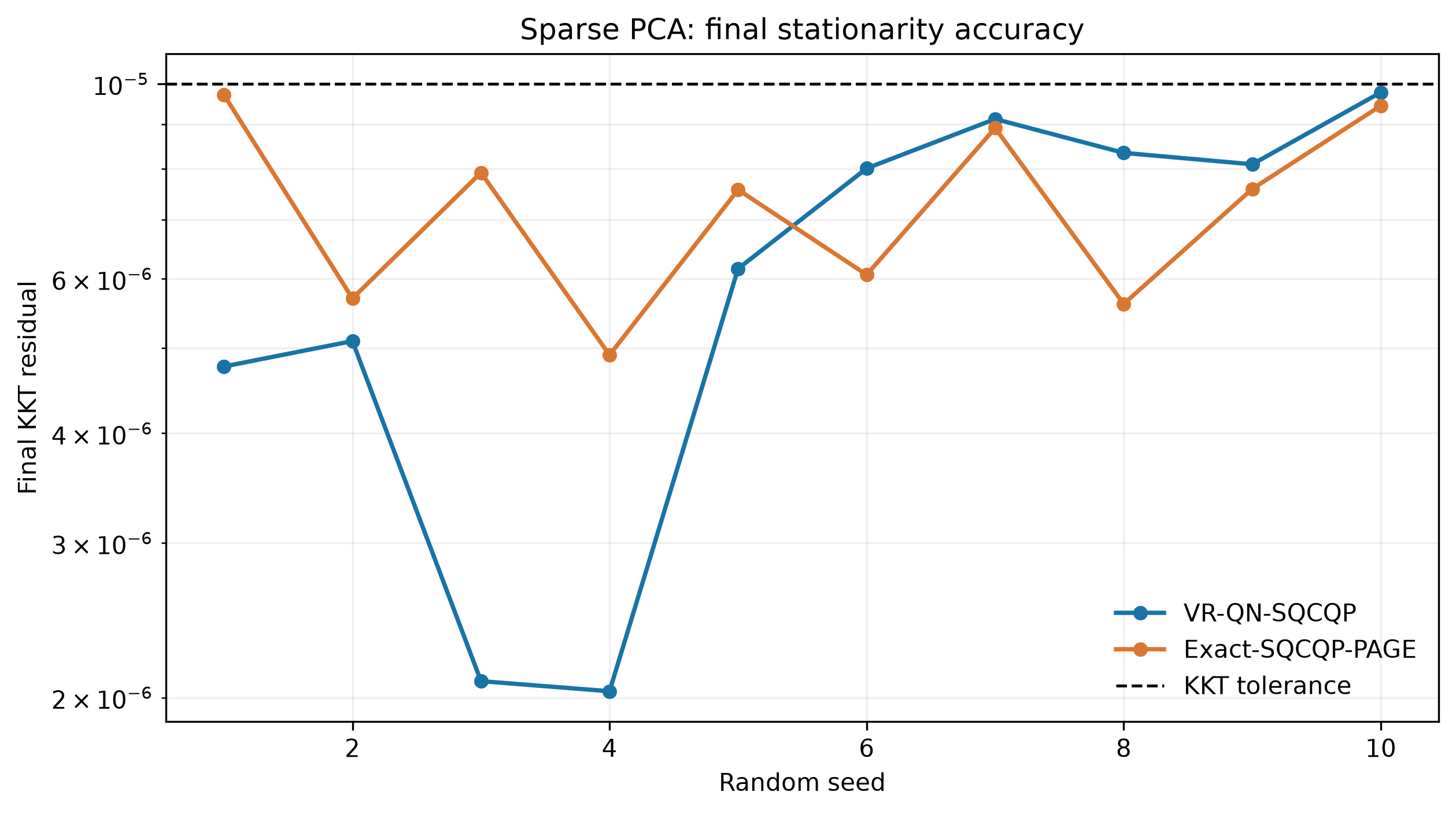}{\textbf{(e)} Final KKT residuals.}\hfill
\paperfigpanel{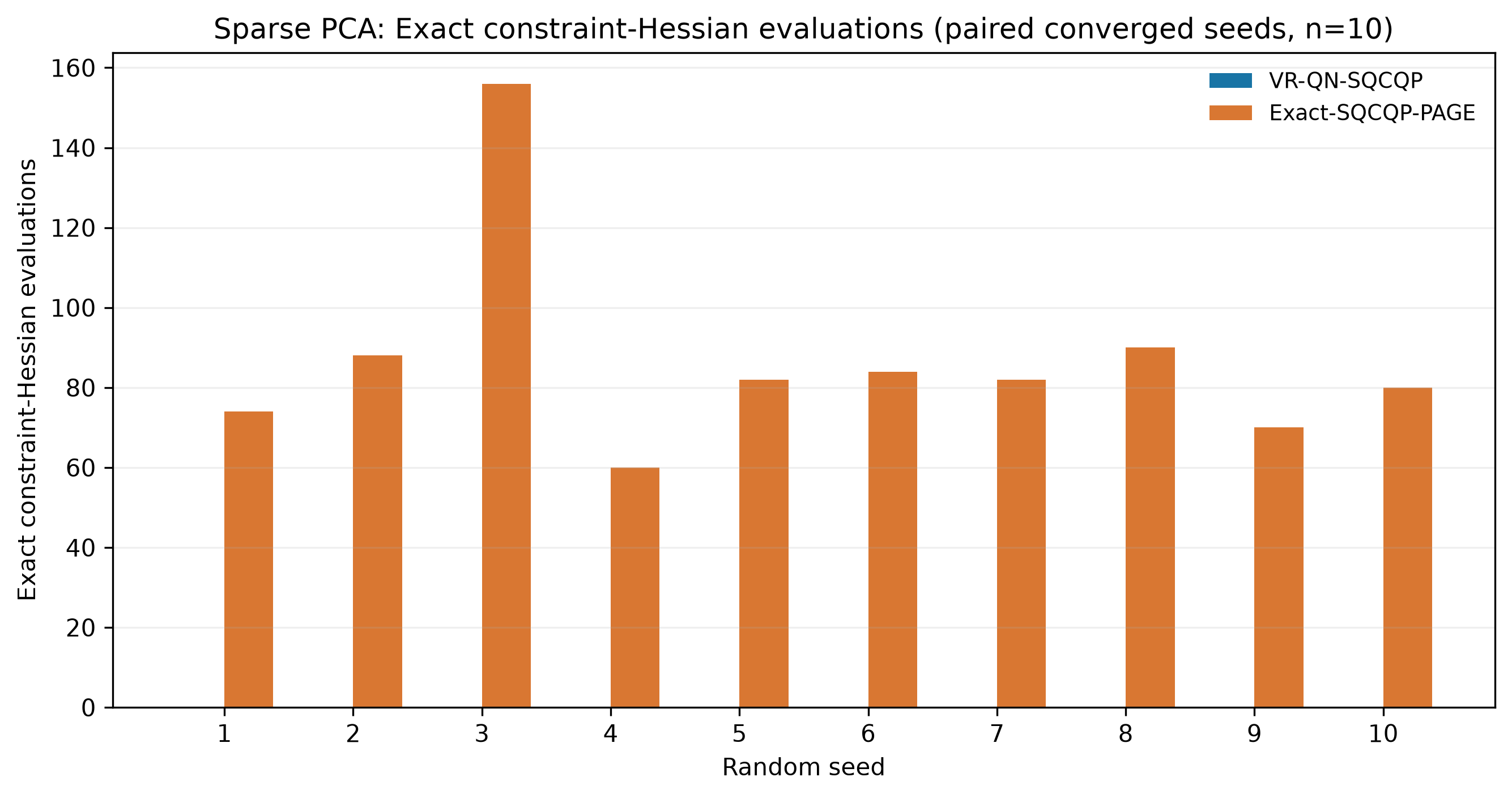}{\textbf{(f)} Exact constraint-Hessian evaluations.}

\par\vspace{3mm}

\paperfigpanel{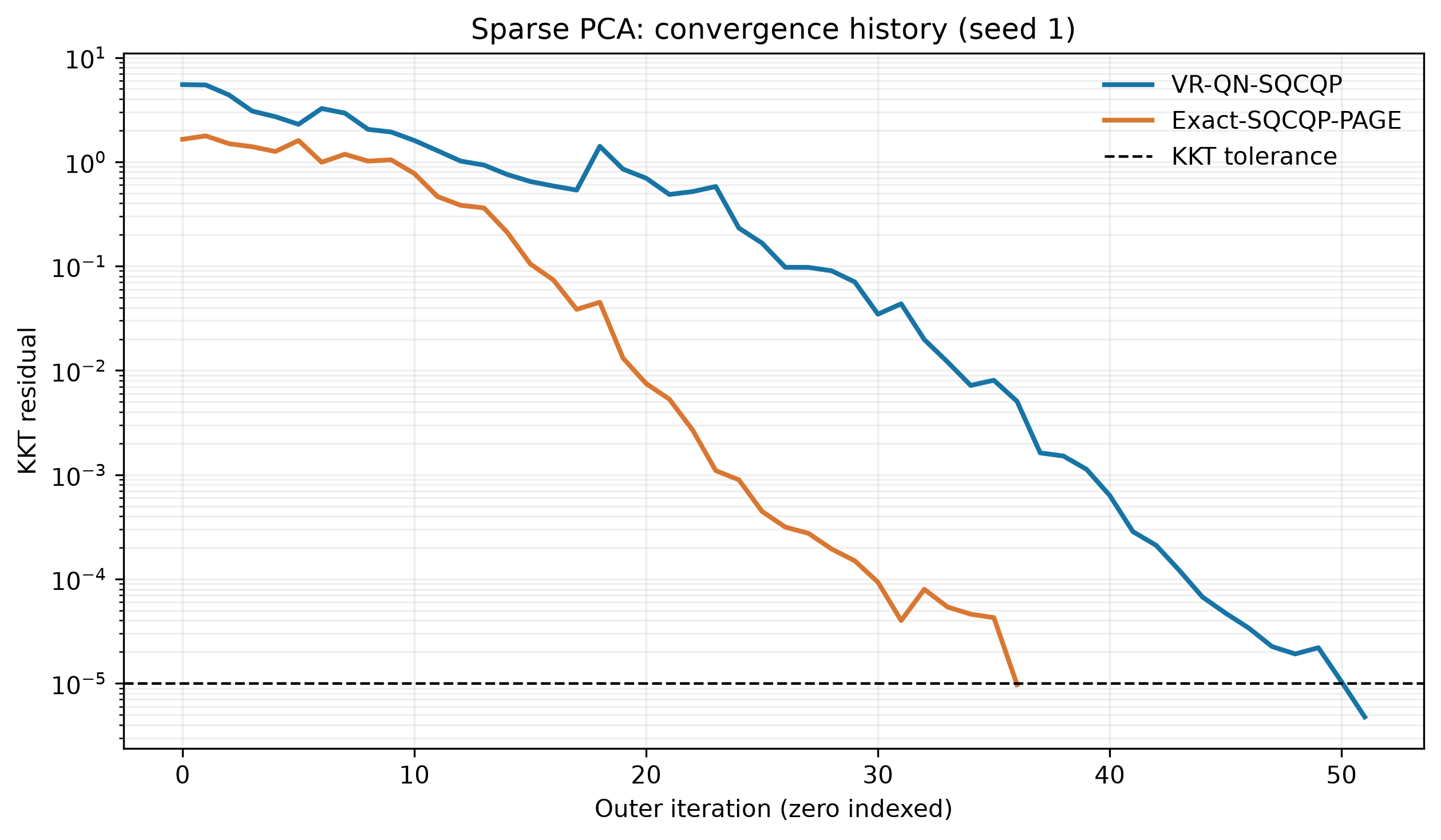}{\textbf{(g)} KKT convergence history.}\hfill
\paperfigpanel{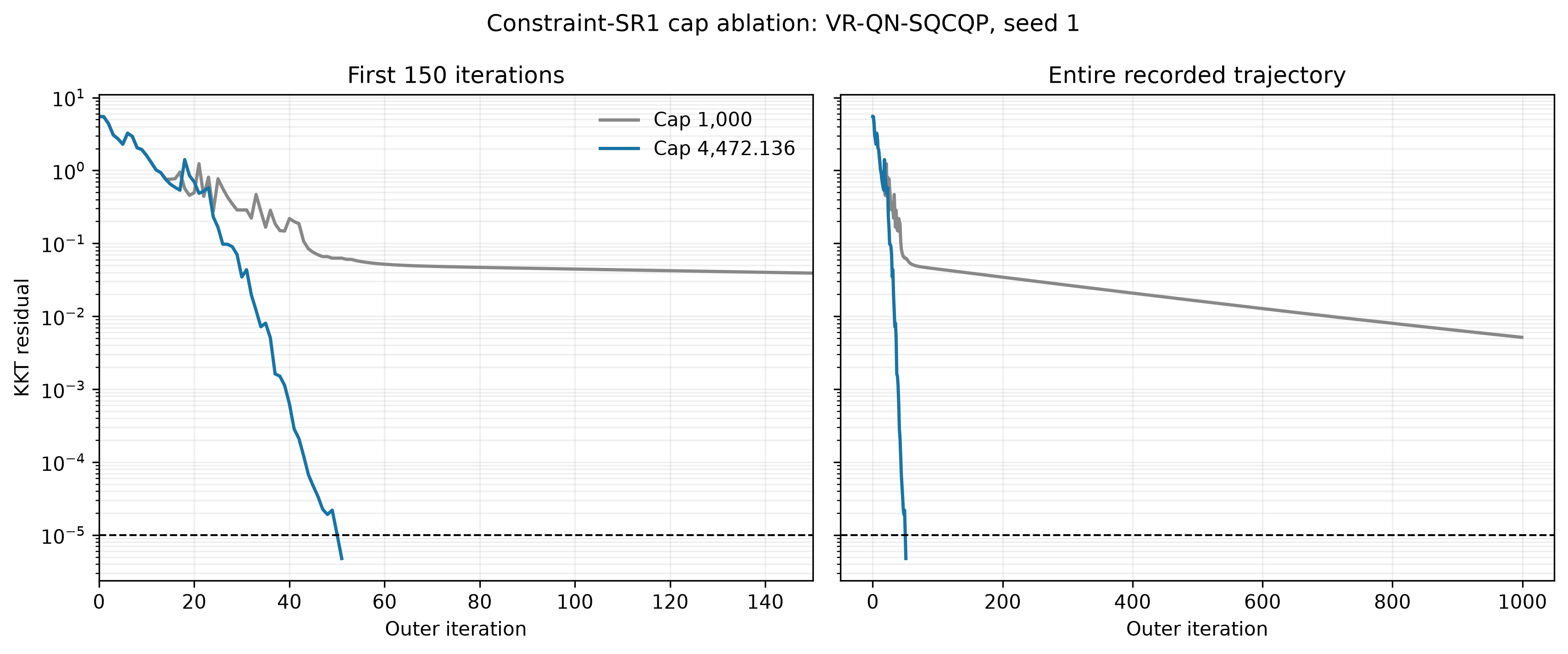}{\textbf{(h)} Constraint-SR1 safeguard ablation.}

\caption{
Numerical performance on full-data Sparse PCA
with $N=253\,680$ observations and $n=20$
predictors. Eight comparisons are displayed:
(a) convergence reliability,
(b) computational runtime,
(c) cumulative QCQP solves,
(d) gradient-oracle samples,
(e) final KKT residuals,
(f) exact constraint-Hessian evaluations,
(g) KKT convergence history for seed 1,
and (h) constraint-SR1 safeguard sensitivity.
The KKT tolerance is fixed at $10^{-5}$.
The safeguard experiments use caps of 1000
and 4472.136, with maximum iteration budgets
of 1000 and 300, respectively.
}

\label{fig:sparse_results}
\end{figure*}
\clearpage
\begin{figure*}[p]
\centering
\paperfigpanel{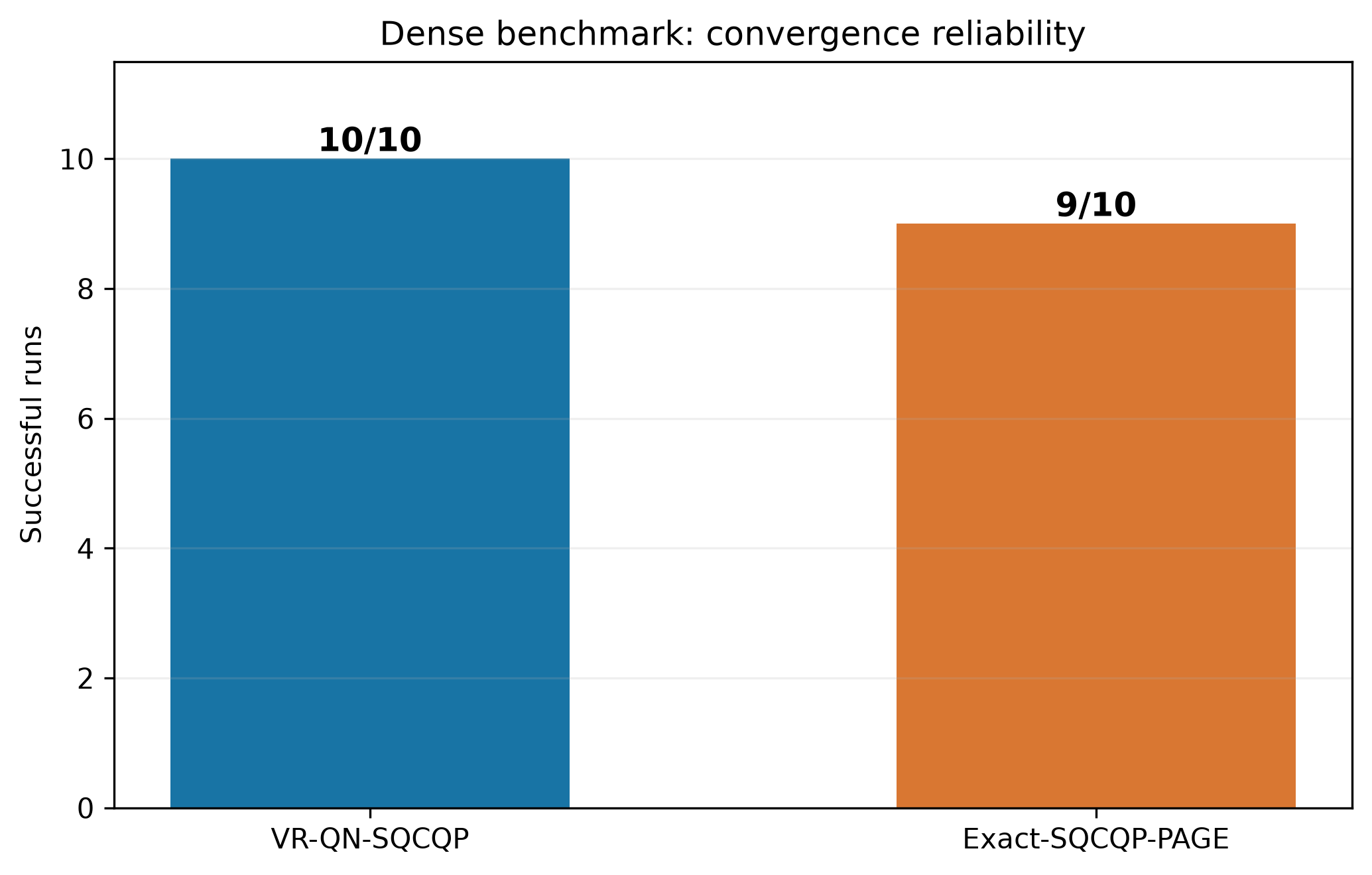}{\textbf{(a)} Convergence reliability.}\hfill
\paperfigpanel{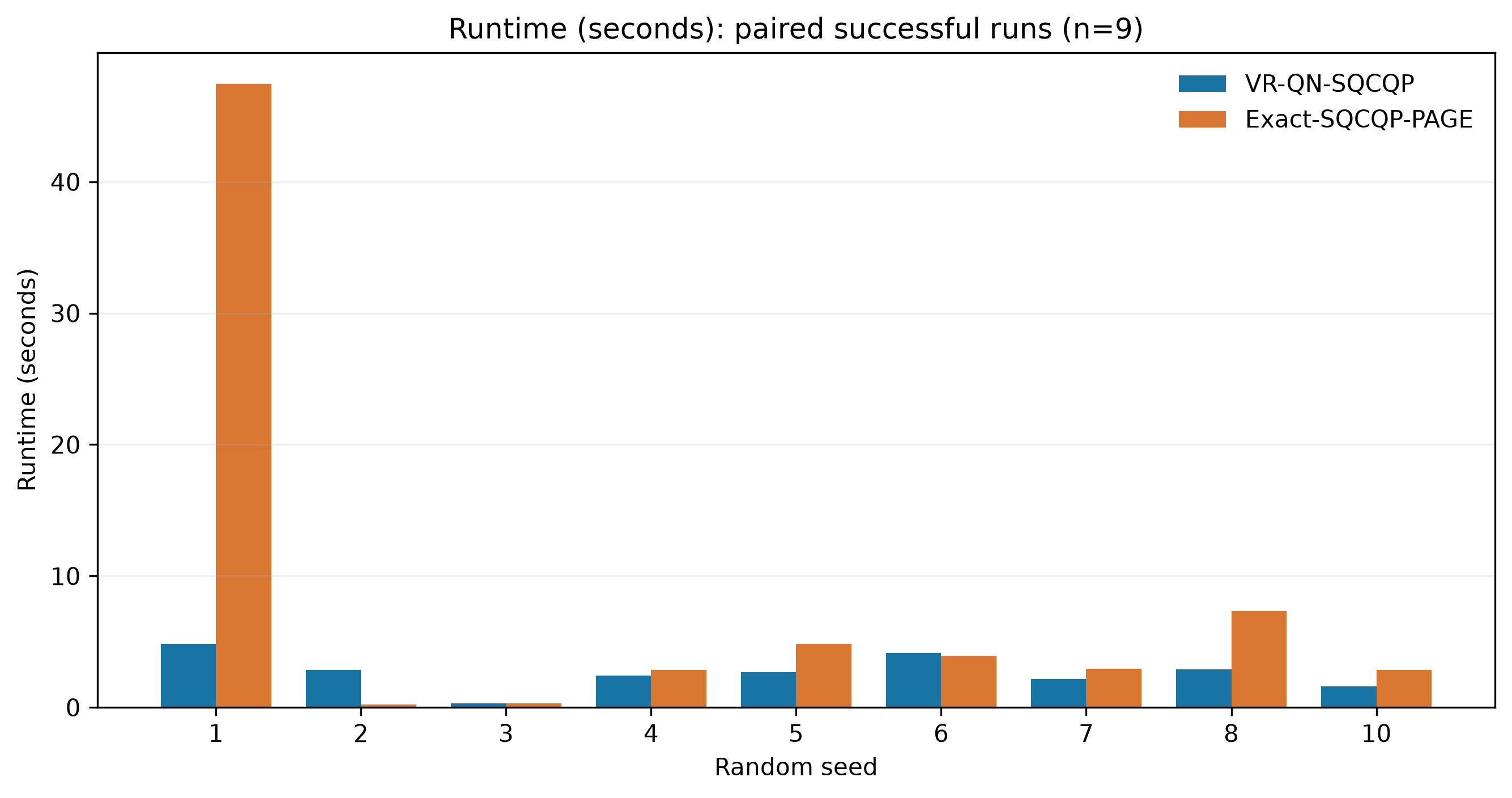}{\textbf{(b)} Computational runtime.}

\par\vspace{3mm}

\paperfigpanel{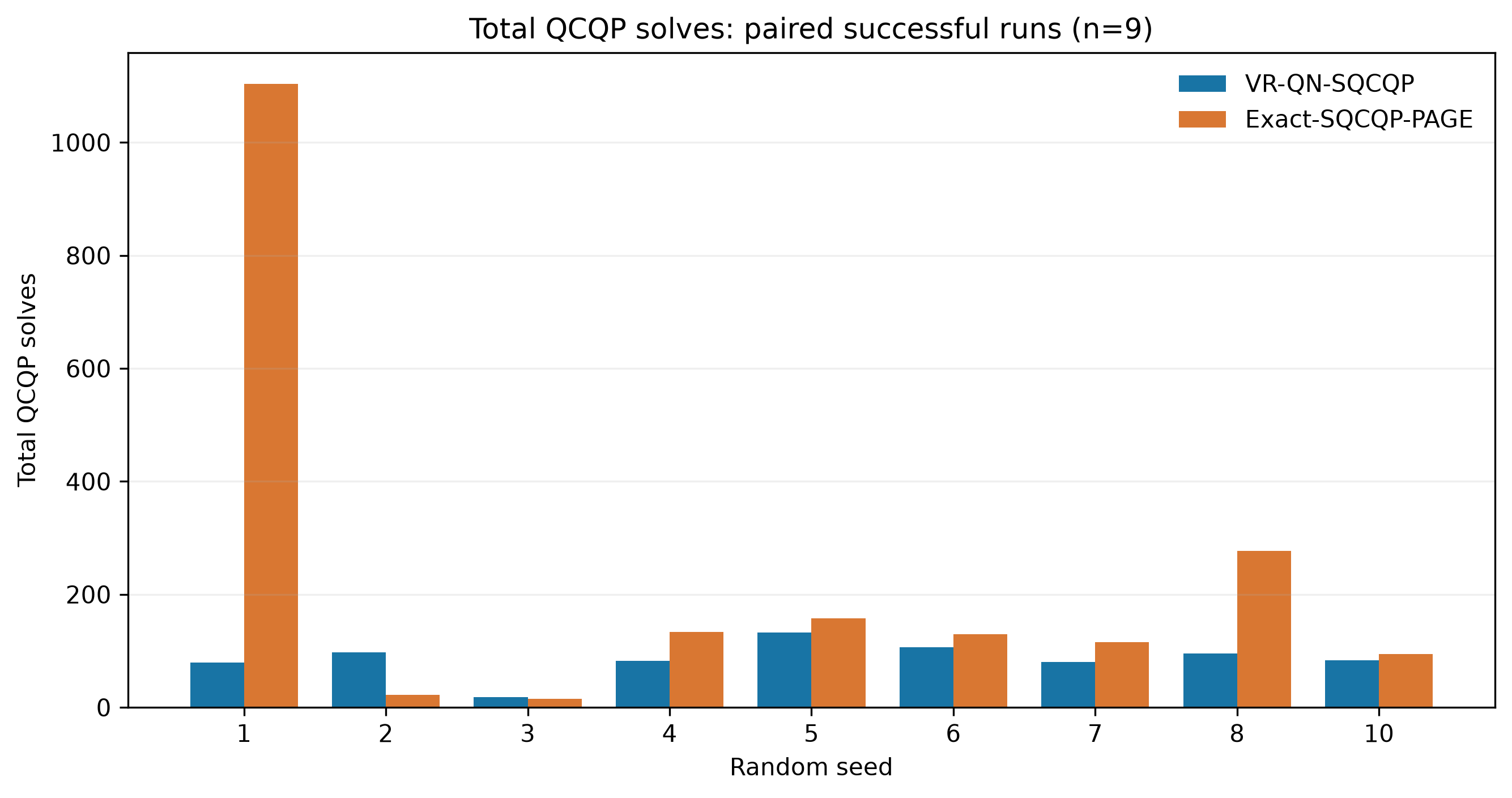}{\textbf{(c)} Total QCQP solves.}
\hfill
\paperfigpanel{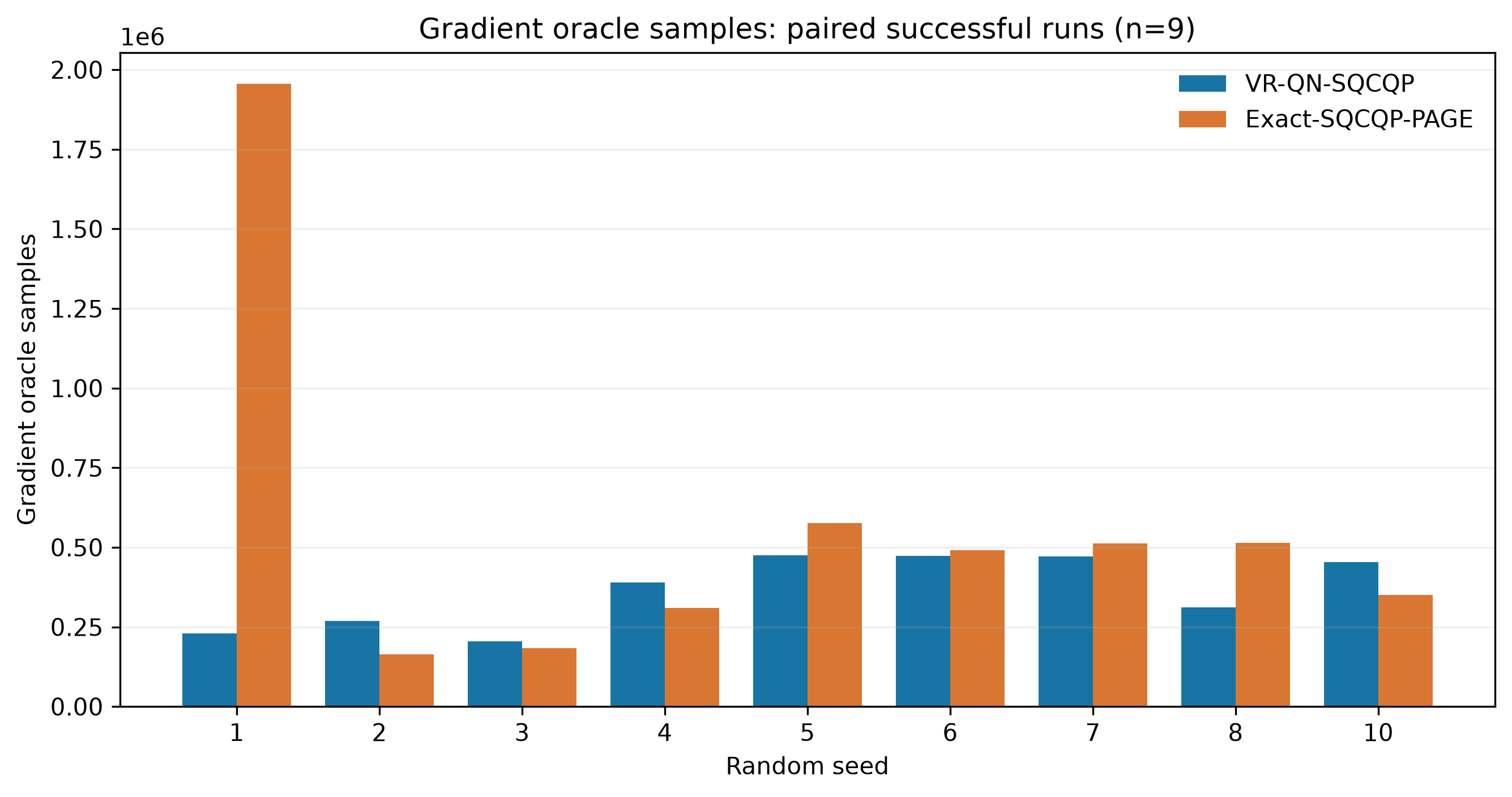}{\textbf{(d)} Gradient-oracle samples.}

\par\vspace{3mm}

\paperfigpanel{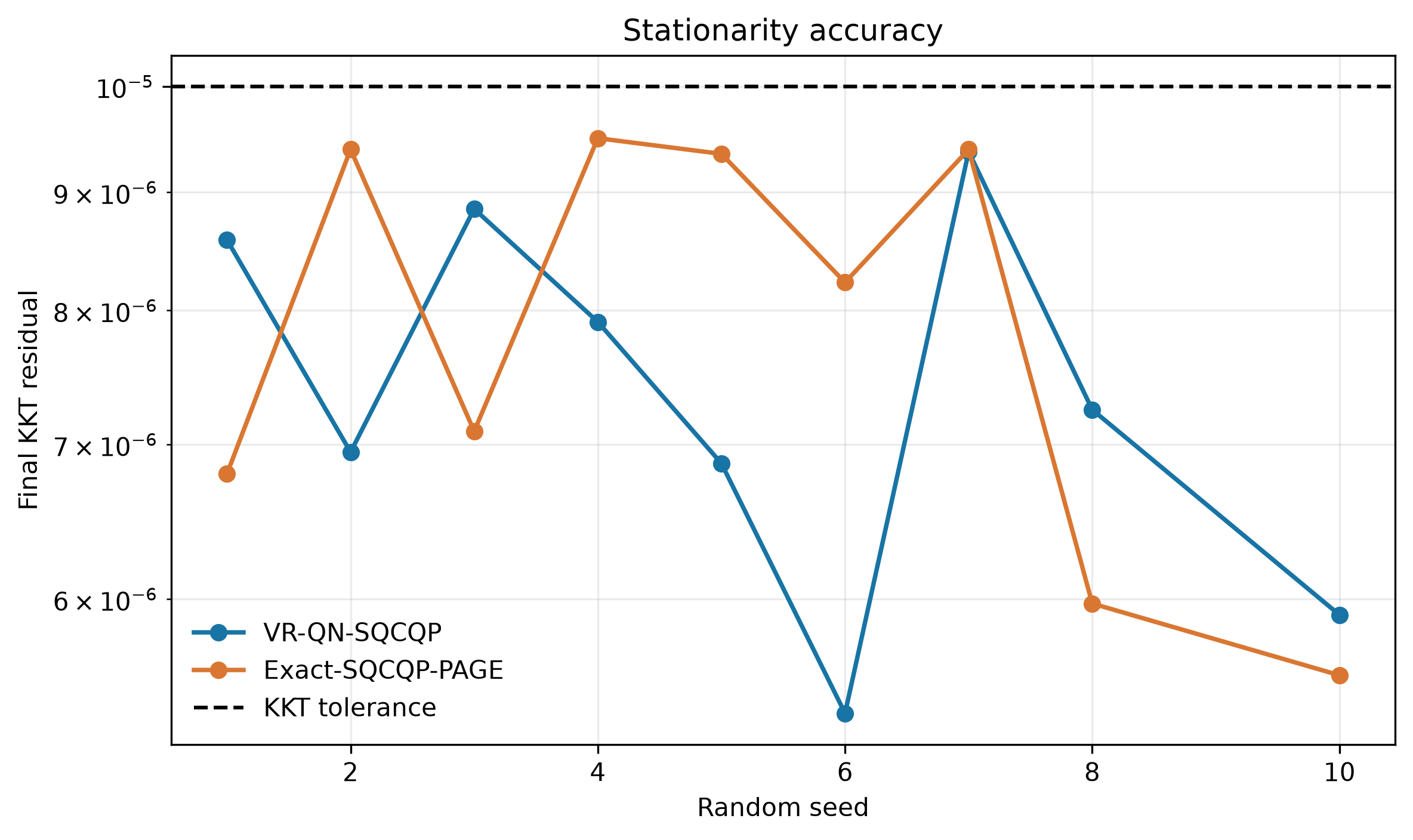}{\textbf{(e)} Final KKT residuals.}
\hfill
\paperfigpanel{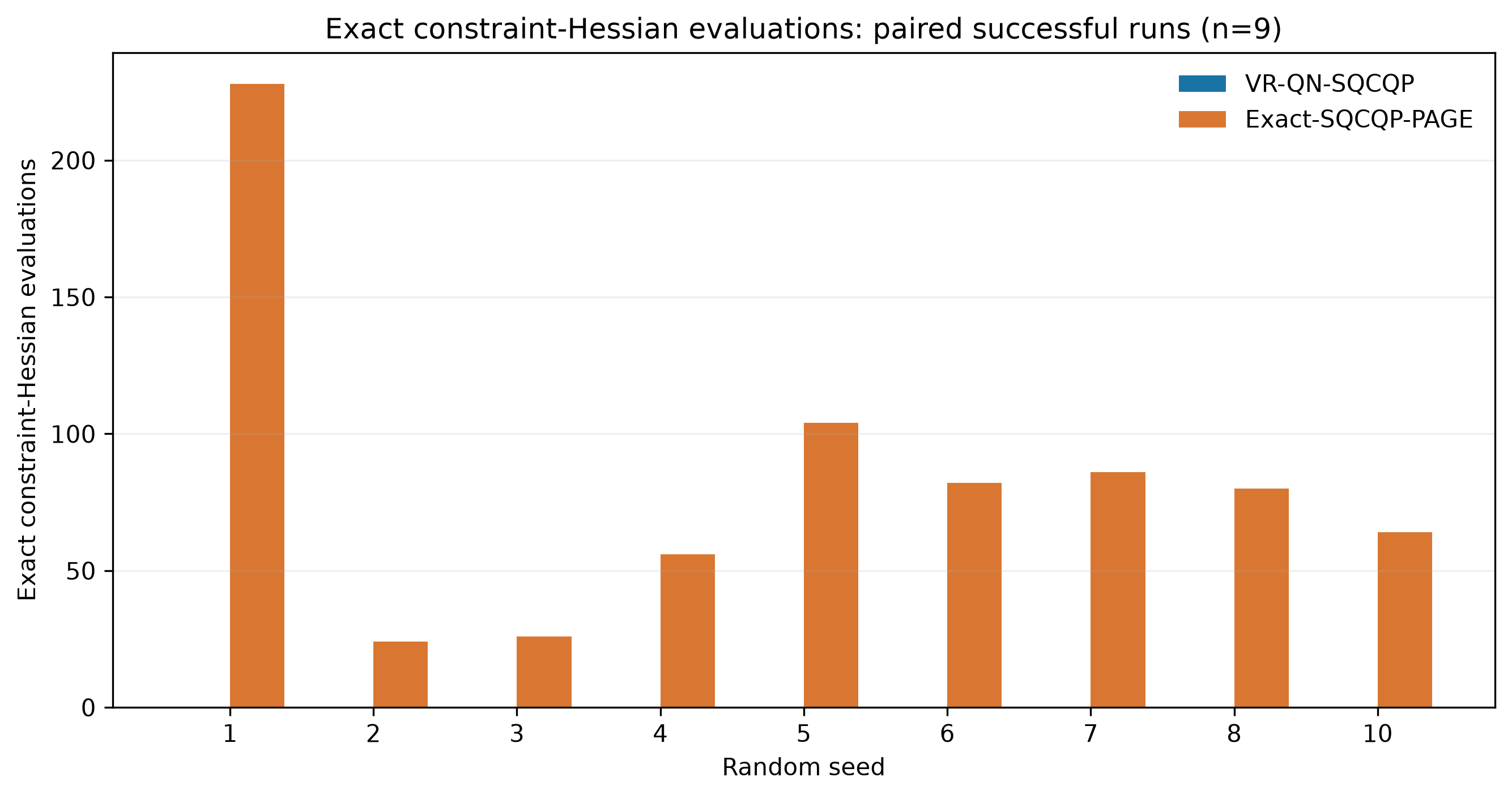}{\textbf{(f)} Exact constraint-Hessian evaluations.}

\par\vspace{3mm}

\paperfigpanel{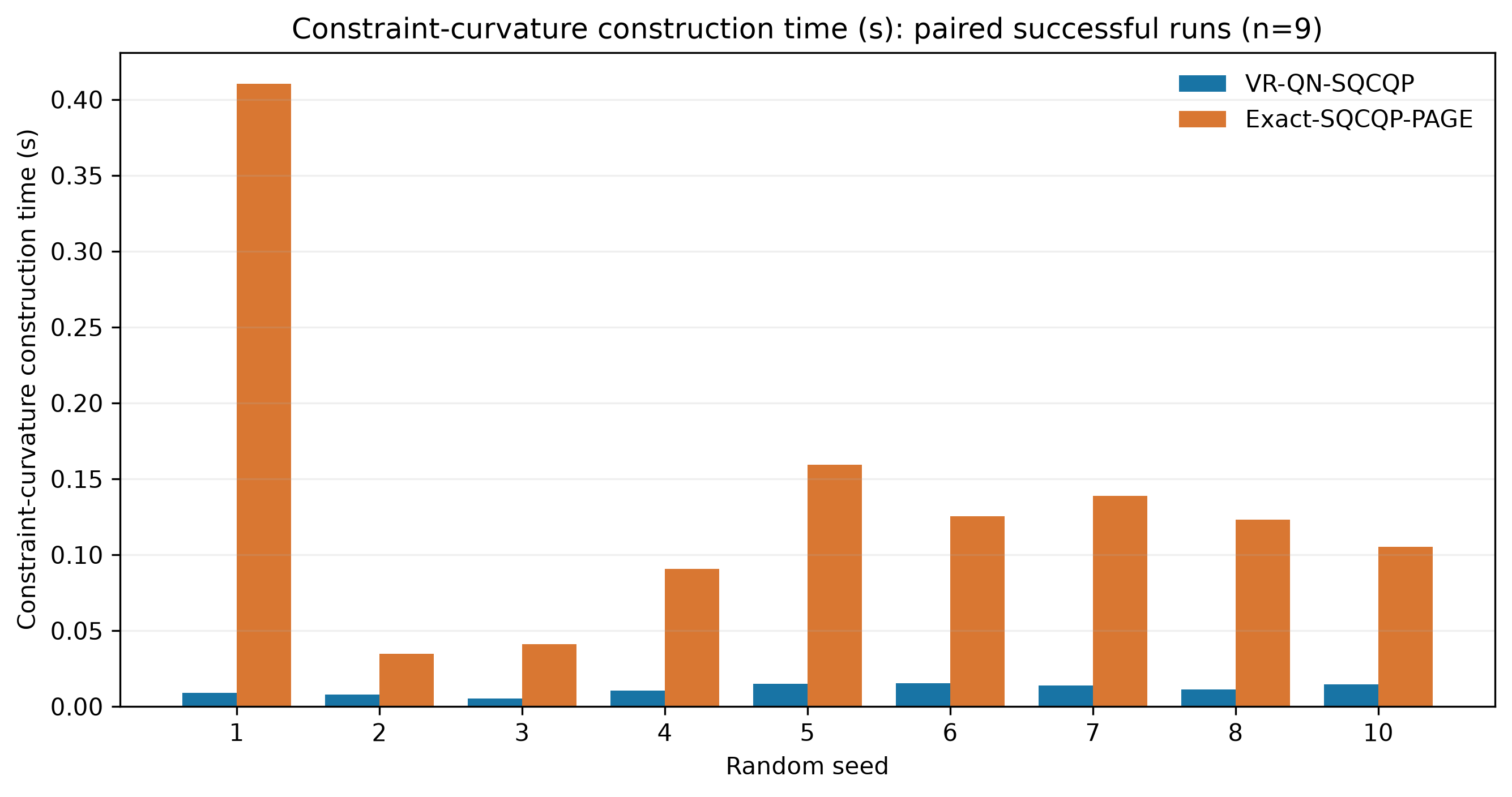}{\textbf{(g)} Constraint-curvature construction time.}\hfill
\paperfigpanel{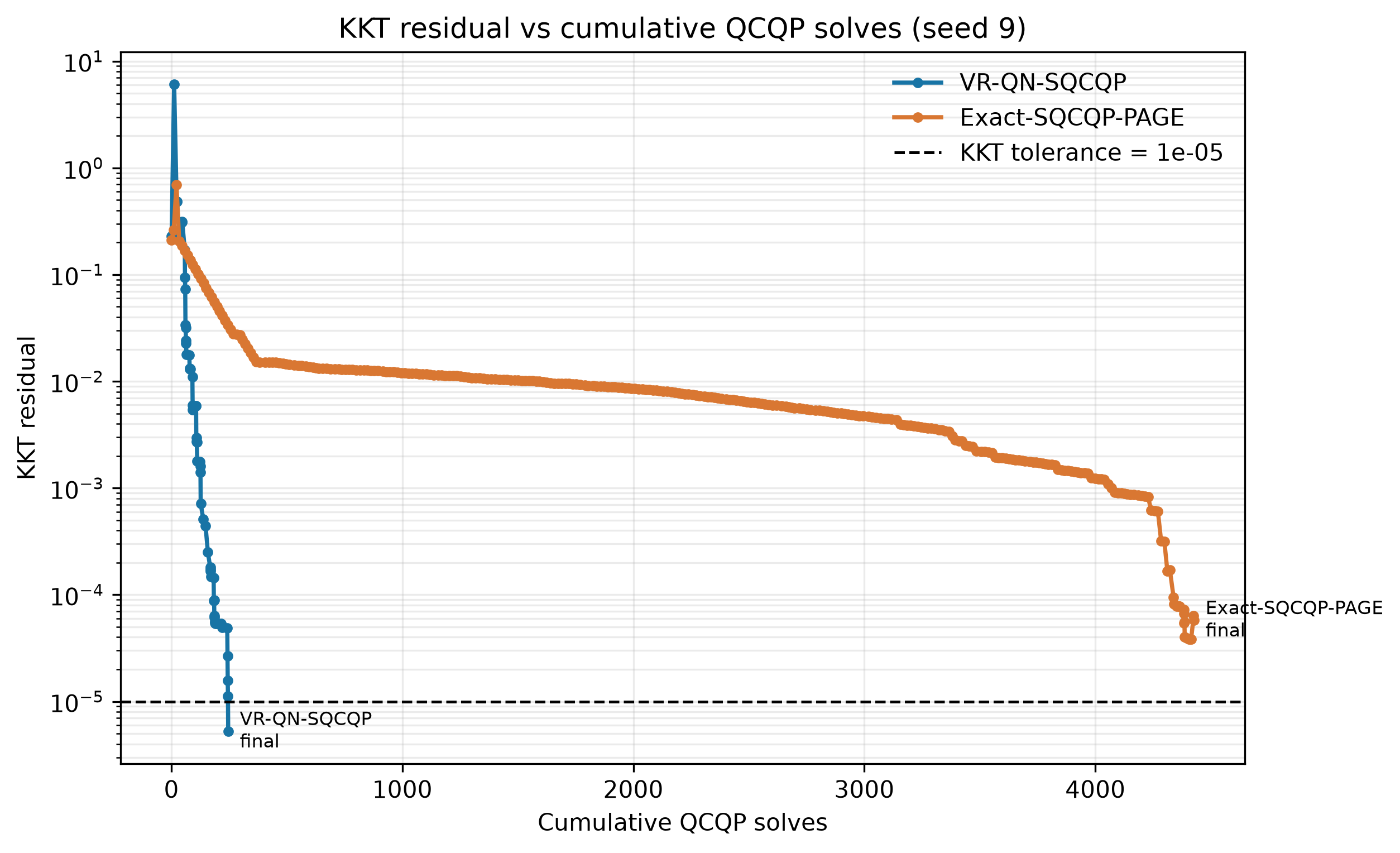}
{\textbf{(h)} KKT residual versus QCQP solves.}

\caption{
Numerical performance on the dense nonlinear
benchmark with $N=20\,000$, $n=100$, and
$q=2000$. The panels illustrate
(a) convergence reliability,
(b) computational runtime,
(c) cumulative QCQP solves,
(d) gradient-oracle samples,
(e) final KKT residuals,
(f) exact constraint-Hessian evaluations,
(g) cumulative constraint-curvature construction
time, and (h) KKT residual versus cumulative
QCQP solves for seed 9.
Ten random seeds are considered.
Computational comparisons are performed
over the nine jointly converged seeds,
whereas convergence reliability is evaluated
over all ten seeds.}

\label{fig:dense_results}
\end{figure*}
\clearpage
\section{Conclusion and Future Scope}\label{sec:conclusion}
A variance-reduced stochastic quasi-Newton SQCQP
(VR-QN-SQCQP) method has been developed for nonlinear constrained stochastic optimization without requiring exact constraint-Hessian evaluations. By incorporating PAGE gradient estimation, safeguarded symmetric-indefinite SR1 curvature approximations, and merit-based globalization, quadratic constraint models have been constructed using only first-order constraint information. Conditional global convergence and local second-order consistency properties have been established under suitable regularity assumptions. Numerical experiments on Sparse PCA and dense nonlinear constrained problems have demonstrated that comparable solution accuracy can be achieved while eliminating exact constraint-Hessian evaluations. Substantial reductions in curvature-construction cost have been observed for dense nonlinear constraints, although overall computational efficiency remains problem dependent. As future research directions, adaptive curvature-safeguarding strategies and more efficient nonconvex QCQP solvers may be investigated to reduce computational overhead. Extensions to stochastic nonlinear constraints, higher-dimensional problems, and stronger local convergence guarantees may also be explored to further enhance the applicability of the proposed framework.

\section*{Declarations}
\begin{itemize}
    \item {\bf Funding:} This article is partially funded by the NBHM research project entitled ``{\it Higher Order Stochastic Methods for Large Scale Multi-objective Optimization Problems}''. Author {\it Md Abu Talhamainuddin Ansary} is thankful to NBHM for this funding.
    \item  {\it Conflict of interest/Competing interests:} The authors declare no conflict of interest.
\item {\it Ethics approval and consent to participate:} Not applicable
\item {\it Consent for publication:} Not applicable
\item {\it Data availability:} The dataset used in this study was obtained from Kaggle and is publicly available at: https://www.kaggle.com/datasets/alexteboul/diabetes-health-indicators-dataset.
\item {\it Materials availability:} Not applicable
\item {\it Code availability:} Codes will be made available from corresponding author based on requests.
\item {\it Author contribution:} Each author has contributed equally 
\end{itemize}

\end{document}